\documentclass[11pt]{amsart}
\usepackage[dvips]{graphicx}
\usepackage{amsmath,graphics}
\usepackage{amsfonts,amssymb}
\theoremstyle{plain}
\newtheorem*{theorem*}{Theorem}
\newtheorem*{lemma*} {Lemma}
\newtheorem*{corollary*} {Corollary}
\newtheorem*{proposition*}{Proposition}
\newtheorem*{conjecture*}{Conjecture}
\newtheorem{theorem}{Theorem}[section]
\newtheorem{lemma}[theorem]{Lemma}
\newtheorem{corollary}[theorem]{Corollary}
\newtheorem{proposition}[theorem]{Proposition}
\newtheorem{conjecture}[theorem]{Conjecture}

\theoremstyle{remark}

\newtheorem*{remark}{Remark}

\theoremstyle{definition}
\newtheorem{defn}[theorem]{Definition}

\def\gl{\mbox{GL}}  \def\F{\Bbb{F}} \def\Z{\Bbb{Z}} \def\R{\Bbb{R}} 
\def\N{\Bbb{N}}  \def\l{\lambda}  
 \def\a{\alpha} \def\g{\gamma} \def\tor{\mbox{Tor}} \def\bp{\begin{pmatrix}}
\def\sm{\setminus} \def\ep{\end{pmatrix}} \def\bn{\begin{enumerate}} \def\Hom{\mbox{Hom}}
 \def\rank{\mbox{rank}} \def\div{\mbox{div}} \def\en{\end{enumerate}}
\def\ba{\begin{array}} \def\ea{\end{array}}  
   \def\a{\alpha} \def\b{\beta} \def\ti{\tilde}
  \def\im{\mbox{Im}} 
  
\def\ker{\mbox{Ker}}\def\be{\begin{equation}} \def\ee{\end{equation}} 
  \def\ord{\mbox{ord}} 
 \def\hom{\mbox{Hom}}  \def\gcd{\mbox{gcd}}
 \def\aut{\mbox{Aut}}  
   \def\zgt{\Z[G][t^{\pm 1}]}

\def\twialexg{\Delta^{\a}_{\phi}(t)}

\def\zt{\Z[t^{\pm 1}]}    
\def\w{\omega}   
    \def\fr12{\frac{1}{2}} \def\z12{\Z[\fr12]} 
 \def\rf{R[F]} 

\def\fpt{\F_p[t^{\pm 1}]}
\def\fp{\F_p}

\def\tpm {[t^{\pm 1}]}

\begin{document}

\title{Twisted Alexander polynomials and symplectic structures}
\author{Stefan Friedl}
\address{Universit\'e du Qu\'ebec \`a Montr\'eal, Montr\'eal, Qu\`ebec}
\email{friedl@alumni.brandeis.edu}
\author{Stefano Vidussi}
\address{Department of Mathematics, University of California,
Riverside, CA 92521, USA} \email{svidussi@math.ucr.edu} \thanks{S. Vidussi was partially
supported by NSF grant \#0629956.}
\date{August 16, 2007}
\begin{abstract} Let $N$  be a closed, oriented $3$--manifold.
A folklore conjecture states that $S^{1} \times N$ admits a symplectic structure
only if $N$ admits a fibration over the circle. The purpose of this paper is to provide evidence
to this conjecture studying suitable twisted Alexander polynomials of $N$, and showing that
their behavior is the same as of those of fibered $3$--manifolds. In particular, we will obtain
new obstructions to the existence of symplectic structures and to the existence of symplectic
forms representing certain cohomology classes of $S^{1} \times N$. As an application of these
results we will show that $S^{1} \times N(P)$ does not admit a symplectic structure, where
$N(P)$ is the $0$--surgery along the pretzel knot $P = (5,-3,5)$, answering a question of Peter
Kronheimer.
\end{abstract}
\dedicatory{Dedicated to the memory of Jerry Levine}
\maketitle

\section{Introduction}

Let $N$ be a 3--manifold. Throughout the paper, unless otherwise stated, we will assume that all
3--manifolds are closed, oriented and connected.  A 4--manifold is called \textit{symplectic} if it admits a closed,
non--degenerate 2--form $\w$ (cf. \cite{GS99}).

Thurston \cite{Th76} showed  that if $N$ admits a fibration over $S^{1}$, then $S^{1} \times N$ admits a symplectic structure.
It is natural to ask whether the converse of this statement holds true; interest on this question was motivated by Taubes' observations on this problem, combined with his results on the Seiberg-Witten invariants of symplectic $4$-manifolds (see \cite{Ta98}). We can state this problem in the following form:

\begin{conjecture} \label{conjfolk}
Let $N$ be a 3--manifold. If $S^1\times N$ is symplectic, then there exists $\phi\in H^1(N)$
such that $(N,\phi)$ fibers over $S^1$.
\end{conjecture}

Here we say that $(N,\phi)$ fibers over $S^{1}$ if the homotopy class of maps $N \to S^1$
determined by $\phi \in H^{1}(N) = [N,S^{1}]$ contains a representative that is a fiber bundle
over $S^{1}$; in that case, we will also say that $\phi$ is a \textit{fibered class}.
The purpose of this paper is to use twisted Alexander polynomials to investigate
Conjecture \ref{conjfolk}.

The strategy of this paper is to show that, if $S^{1} \times N$ is symplectic, then $N$
satisfies properties that hold for fibered $3$--manifold, obtaining this way evidence for
Conjecture \ref{conjfolk} and obstructions to the existence of symplectic forms.
This approach is made possible by the fact that we can translate Taubes' constraints on the Seiberg-Witten invariants
of a symplectic $4$-manifold of the form $S^{1} \times N$ (hence of $N$) in terms of the properties of the Alexander polynomial of $N$, using the relation between Seiberg-Witten and Alexander invariants of a $3$--manifold determined by Meng and Taubes in \cite{MeTa}.
This strategy
was first used by Kronheimer who proved  in \cite{Kr99} that if the $0$--surgery $N(K)$ along a
knot $K \subset S^{3}$ is such that $S^{1} \times N(K)$ is symplectic, then the Alexander
polynomial of $K$ is monic and its degree coincides with twice the genus of $K$. These
conditions are known to hold for fibered knots. This approach was generalized by the second
author in \cite{Vi03} to all irreducible $3$--manifold (see Section \ref{intro} for
details). This paper uses the twisted Alexander polynomial associated to an epimorphism of the
fundamental group of a $3$--manifold onto a finite group to strengthen these conclusions. This
way we obtain further evidence for Conjecture \ref{conjfolk} and new obstructions to the
existence of symplectic forms.

The main result in this sense is Theorem \ref{strongsymp}. For sake of exposition, we quote here
a result that is a slightly weaker version of Theorem \ref{strongsymp}:

\begin{proposition*} Let $N$ be an irreducible $3$--manifold such that $S^{1} \times N$ admits a symplectic
structure. Then there exists a non--trivial  $\phi \in H^{1}(N)$ such that for any epimorphism
$\alpha: \pi_{1}(N) \rightarrow G$  onto a finite group $G$ the associated $1$--variable twisted
Alexander polynomial $\Delta_{N,\phi}^{\a}\in \zt $  is monic and
\[ \deg\,\Delta_{N,\phi}^{\a} = |G|\, \|\phi\|_{T} + 2 \div \, \phi_{G}.\]
Here we denote by $\phi_{G}$ the restriction of $\phi:\pi_1(N)\to \Z$ to $\ker\{\a: \pi_1(N)\to
G\}$, $\div\,  \phi_G\in \N$
denotes its divisibility  and $\|\phi\|_{T}$ is the
Thurston norm of $\phi$.
\end{proposition*}

In fact for any 3--manifold we have $\deg\,\Delta_{N,\phi}^{\a} \leq |G| \, \|\phi\|_{T}
+ 2 \, \div \, \phi_{G}$. This inequality generalizes McMullen's inequality for the ordinary
Alexander polynomial (Lemma \ref{normineq}). Note that whenever $\phi$ is a fibered class the
equality holds and $\Delta_{N,\phi}^\a$ is monic (Theorem \ref{mainthmg}). In this sense, this
proposition supports Conjecture \ref{conjfolk}.

In practical terms, the previous proposition can be used to show that certain manifolds of the
form $S^{1} \times N$, where $N$ is a non--fibered $3$--manifold for which the ordinary
Alexander polynomial does not provide an obstruction, do not admit any symplectic form. An
interesting example in this sense is the manifold $N(P)$ given by  $0$--surgery of $S^3$ along the pretzel knot $P
= (5,-3,5)$. This knot (and so $N(P)$) is not fibered, has genus $1$, and its Alexander
polynomial equals $\Delta = t^{2} -3t + 1$, so Kronheimer's results  do not exclude that $S^{1}
\times N(P)$ is symplectic. Kronheimer asked in \cite{Kr99} to determine whether $S^{1} \times
N(P)$ is symplectic or not. We have the following.

\begin{proposition*} Let $P$ be the pretzel knot $(5,-3,5)$,
and let $N(P)$ be the $0$--surgery of $S^3$ along $P$. The manifold $S^{1} \times N(P)$ does not
admit a symplectic structure.
\end{proposition*}

The proof follows by applying the previous proposition for a suitable choice of the epimorphism
$\alpha$. Several other examples, similar to this, are presented. In particular we completely
determine for which knots with up to 12 crossings the manifold  $S^{1} \times N(K)$ admits a
symplectic structure. For many of these knots the ordinary Alexander polynomial does not contain
sufficient information. This leads to verifying the following.
\begin{proposition*}
Conjecture \ref{conjfolk} holds true for a $3$--manifold obtained by $0$--surgery along a knot with up
to $12$ crossings. \end{proposition*}

For fibered manifolds $N$ with $b_{1}(N) > 1$ the first proposition
can be seen as determining new constraints on the size of the symplectic cone of $S^{1} \times
N$, as the K\"unneth component in $H^{1}(N)$ of an integral symplectic form must satisfy the
conclusion of the theorem.

In view of the strength of the constraints on fiberability given by the twisted Alexander
polynomials, we will present (perhaps optimistically) a conjecture (Conjecture \ref{conjcha}),
stating that the constraints associated to all epimorphisms of the fundamental group of a
$3$--manifold onto finite groups determine whether the manifold is fibered. We will prove that
this conjecture implies Conjecture \ref{conjfolk}, and discuss some arguments in favor of the
conjecture.

We finish this introduction with two remarks. First, in this paper we have followed the approach
of formulating  definitions and results in terms of (twisted) Alexander polynomials and the
corresponding Alexander norms. This is convenient because of their wide use, but we would like
to point out that the most natural approach would be to use instead the corresponding
Reidemeister torsion. Several results (particularly for the $b_{1}(N) = 1$ case) have a more
elegant presentation when formulated in terms of torsion. Second, some of the results concerning
the fiberability of a $3$--manifold that are presented here appear (although in different form
or generality) in two papers of the first author and Taehee Kim (see \cite{FK05} and
\cite{FK06}) as well as, for the case of knots, in previous work of Cha (see \cite{Ch03}) and
Goda, Kitano and Morifuji (see \cite{GKM05}). However, we have decided for their inclusion both
in light of Conjectures \ref{conjfolk} and \ref{conjcha} (as they specify the actual evidence
to it) and because the proofs presented are different in nature and have often emerged
independently. In particular, the approach presented here is naturally influenced by the
viewpoint of Seiberg-Witten theory, and we would like to point out that all our constraints on
fiberability (and the McMullen-type inequality) can be obtained with the use of Seiberg-Witten
theory alone. In this sense we hope to have progressed (although very modestly) in the direction
of better understanding the ``intriguing circle of ideas" referred to in the introduction of
\cite{MT99}.
\\

\noindent {\bf Organization of the paper:} Section \ref{intro} summarizes what is known on the
relation between Conjecture \ref{conjfolk} and the ordinary Alexander polynomial. In Section
\ref{twistal} we define the twisted Alexander polynomials and discuss some of their properties,
with particular attention to the case of polynomials associated with epimorphisms onto finite
groups. In Section \ref{fini} we determine the new obstructions to the existence of symplectic
structures in terms of twisted Alexander polynomials, and we apply this result to study a family
of  examples in Section \ref{sectionapp}. In Section \ref{section:finitefields} the obstructions
are reformulated in terms of certain twisted Alexander polynomials with coefficients in a finite
field, which are usually easier to compute. This allows us, in Section \ref{section:ex}, to
decide Kronheimer's example and determine in the affirmative Conjecture \ref{conjfolk} for knots with up to 12 crossings.
Sections \ref{sectionconj} contains some remarks on our conjecture.
\\

\noindent {\bf Notations and conventions:} We assume that all manifolds are closed, oriented and
connected. All homology groups and cohomology groups are with respect to $\Z$--coefficients,
unless specifically said otherwise. For a knot $K$ in $S^3$ we denote the result of zero
framed surgery along $K$ by $N(K)$. For a submanifold $X$ of a manifold $M$ we write $\nu X$ for an
open tubular neighborhood of $X$ in $M$. Given a prime $p$ the finite field with $p$ elements is
denoted by $\F_p$. For a 3--manifold $N$ we write $H = H_1(N)/\tor(H_1(N))$, i.e. $H$ is the
maximal free abelian quotient of $H_1(N)$. Given a free abelian group $F$ (by which we always
mean a finitely generated free abelian group) we use the canonical isomorphisms to identify
$H^1(N;F) = \Hom(\pi_1(N),F) = \Hom(H_1(N), F)=\hom(H,F)$. Given a domain $R$ we denote its
quotient field by $Q(R)$. Given a (co)homology class $a$ on $N$ we denote its Poincar\'e dual by
$PD_Na$.
\\

\noindent {\bf Acknowledgment:} The authors would like to thank John Hempel, Taehee Kim and
Jerry Levine for helpful discussions.

\section{Symplectic products and the ordinary Alexander polynomial}
\label{intro}

Our first goal is to discuss some properties satisfied by fibered $3$--manifolds. In order to
present these properties, we will begin by reviewing some definitions that are relevant in what
follows.

First, the Thurston (semi)norm (cf. \cite{Th86}):  Let $\phi \in H^{1}(N)$.
The \emph{Thurston norm} of $\phi$ is defined as
\[
||\phi||_{T}=\min \{ -\chi(\hat{S})\, | \, S \subset M \mbox{ properly embedded surface dual to
}\phi\}
\] where $\hat{S}$ denotes the result of discarding all connected components of $S$ with positive Euler
characteristic.
Thurston showed that  we  can uniquely extend the
definition to a continuous norm on  $H^{1}(N;\R)$ by requiring linearity. We denote the norm
ball $\{ \phi\in H^1(N;\R) | \, ||\phi||_T\leq 1\}$ by $B_T$.

A specific case has particular relevance for us: Given an oriented non--trivial knot $K \subset
S^{3}$, let $\phi\in H^1(N(K))\cong \Z$ be a generator. By a result of Gabai
\cite[Theorem~8.8]{Ga87b} we have $||\phi||_T=2 g(K)-2$. Here $g(K)$ is the usual genus of the
knot $K$, defined as the minimal genus of a surface bounding $K$.

Second, the Alexander polynomial and the Alexander (semi)norm: We can associate to $N$ the
multivariable polynomial $\Delta_{N} \in \Z[H]$, where $H = H_1(N)/\tor(H_1(N))$, which we can
assume, by \cite[p.~176]{Tu86}, to be symmetric (cf. Section \ref{sectionufd} for details). We
can therefore write $\Delta_N=\sum_{g\in H} a_g g$, where $a_{-g}= a_{g}$. The \emph{Alexander norm}
of $\phi \in H^{1}(N;\R)$ is then defined as
\[ \|\phi\|_{A} := 2\, \mbox{max} \{ \phi(g) | g \in \mbox{supp} \Delta_{N} \}. \]
Note that when $\Delta_N=0$ this is understood to mean that $||\phi||_A=0$.
The norm ball $B_{A}$ is a convex (possibly non--compact) polyhedron,
dual to the Newton polyhedron of $\Delta_{N}$ (cf. \cite{McM02} for details). Whenever useful, we
will specify in the notation the manifold to which the Alexander and the Thurston norms refers to.

In \cite{McM02} McMullen has shown that the Alexander norm of a $3$--manifold $N$ provides a
bound by below to the Thurston norm of a cohomology class $\phi \in H^{1}(N)$. McMullen's
inequality has the form
\[ \| \phi \|_{A}  \leq \| \phi \|_{T} + \left\{ \begin{array}{ll} 0 & \mbox{if $b_{1}(N) >
     1$},  \\ \\ 2 \, \div \, \phi  & \mbox{if $b_{1}(N) =
     1$}. \end{array} \right. \]

When the class $\phi$ represents a fibration, this bound becomes an equality, and $\Delta_{N}$
has peculiar properties, as shown from the following theorem, implicitly or explicitly proven,
with different techniques, in \cite{Du01}, \cite{McM02}, \cite{Vi03}, and \cite{FK06}. Remember
that Thurston \cite{Th86} showed that if a class $\phi \in H^{1}(N)$ is fibered then it lies in
the cone over some top--dimensional face $F_{T}$. (Here and in what follows we will
always consider \textit{open} faces for polyhedra.)
Furthermore a class $\phi \in H^{1}(N)$ is
fibered if and only if all the lattice points in the cone $\R_{+} F_{T}$ over the same
top--dimensional face $F_{T}$ of the Thurston unit ball are fibered. The case of completely
degenerate norm is included in this and later statements by interpreting $F_{T}$ as the unique
``face at infinity", dual to the trivial element of $H$, so that $\R_{+} F_{T} = H^{1}(N;\R)
\setminus \{0\}$. Note that in the case $b_1(N)=1$ the (top--dimensional) faces are two discrete points,
in particular they are open.
 Denote by $e_{F_{T}} \in H$ the (Poincar\'e dual of the) Euler class of the
fibrations in $F_{T}$. (Strictly speaking, this is the image of the Euler class under the
quotient map $H_{1}(N) \rightarrow H$.)

\begin{theorem}\label{thmmcmullen} Let $N \neq S^{1} \times S^{2}$ be a fibered $3$--manifold, and
let $F_{T}$ be a fibered face of the Thurston unit ball $B_{T}$. Then the following holds true:
\bn
\item For all $\phi \in \R_{+}F_{T} \cap H^{1}(N)$ \[ \| \phi \|_{A}  = 2 \phi(\lambda) = \| \phi
\|_{T} + \left\{ \begin{array}{ll} 0 & \mbox{if $b_{1}(N) >
     1$},  \\ \\ 2 \, \div \, \phi  & \mbox{if $b_{1}(N) =
     1$}, \end{array} \right. \]
where $\lambda \in H$
is a vertex of the Newton polyhedron of $\Delta_{N}$ with coefficient $\pm
1$. If   $b_{1}(N) > 1$, then $\lambda = -\frac12 e_{F_{T}}$ and $F_{T} \subset F_{A}$ (where
$F_{A}$ is the top--dimensional face of the Alexander unit ball dual to $\lambda$) and if $b_{1}(N)
= 1$, then $\lambda = -\frac12 e_{F_{T}} + 1$.
\item For all $\phi \in \R_{+}F_{T} \cap H^{1}(N)$,
$\Delta_{N,\phi}$ is monic and $\deg\, \Delta_{N,\phi} = \|\phi\|_{T} + 2 \div \, \phi$.
\en
\end{theorem}

Here $\Delta_{N,\phi}$ is the $1$--variable Alexander polynomial associated to $\phi$, viewed as
element of $Hom(H,\Z)$ (cf. Section \ref{sectionufd} for details). Furthermore note that in the
case $b_1(N)=1$ the choice of a face $F_T$ corresponds to the choice of an isomorphism $H\to
\Z$, i.e. we can identify $H$ and $\Z$. The last equation of part (1) in the preceding theorem
has to be understood using this identification.

These properties (that generalize classical results of Neuwirth for knots) are short of
characterizing fibered classes, but provide useful and often effective constraints for
fiberability.

Let us now assume that $S^{1} \times N$ admits a symplectic form $\w$, whose cohomology class
with no loss of generality can be taken in the integer lattice of $H^{2}(S^{1} \times N;\R)$.
Decompose $[\w]$ according to its K\"unneth decomposition, and let $\varphi$ be its K\"unneth
component in $H^{1}(N)$, normalized so that $[\w] = [dt] \wedge \varphi + \eta$, where $[dt]$ is
the generator of $H^{1}(S^1)$ and $\eta \in H^{2}(N)$. Denote by $K \in H^{2}(S^{1} \times N)$
the canonical class of the symplectic structure: such class is well--known to be the pull-back
of a class in $H^{2}(N)$ that we denote with the same symbol. Evidence to Conjecture
\ref{conjfolk} comes from the fact that such a manifold (under the additional hypothesis of
irreducibility) satisfies the same properties of Theorem \ref{thmmcmullen}, as shown in
\cite{Vi03} (and \cite{Vi99} for $b_{1}(N) = 1$):
\begin{theorem} \label{thmvidussi}
Let $N$ be an irreducible 3--manifold such that $S^1\times N$ admits an integral symplectic
structure with K\"unneth component $\varphi \in H^{1}(N)$. Then there exists a top--dimensional
face $F_{T}$ of the Thurston unit ball with $\varphi \in \R_{+}F_{T}$ for which the same
properties as in (1) and (2) of Theorem \ref{thmmcmullen} hold true, with $\lambda$ equal to
$\frac12 PD_{N} K$ if $b_{1}(N) > 1$ and $\frac12 PD_{N} K + 1$ if $b_{1}(N) = 1$.
\end{theorem}

For future reference and completeness, we summarize the idea of the proof of \cite{Vi03}. This
proof combines (generalizing ideas of Kronheimer, see \cite{Kr98} and \cite{Kr99}) the following
ingredients, connected by the relation between the Seiberg-Witten theory of $S^{1} \times N$ and
the Alexander and Thurston norm of $N$. First, Taubes' results on symplectic $4$--manifolds
impose conditions on the convex hull of basic classes of Seiberg-Witten theory, in particular on
$K$, and its relation with $\varphi$. In particular, $\varphi$ has maximal pairing, among all
elements of $\mbox{supp}\Delta_{N}$, with $\lambda$. This implies that $\varphi \in \R_{+}F_{A}$,
where $F_{A}$ is the face dual to $\lambda$. Second, Donaldson's Theorem guarantees the existence
of a symplectic representative for the class Poincar\'e dual to (a sufficiently high multiple
of) $\w$. Last, we apply to such a representative Kronheimer's refined adjunction inequality for
manifolds of the form $S^{1} \times N$. These constraints prove the equalities of Theorem
\ref{thmvidussi} for the norms of $\varphi$. Using the openness of the symplectic condition we can extend
this result, when $b_{1}(N) \geq 2$, to some open cone in $H^{1}(N;\R)$, determined by symplectic
forms all having canonical class $K$. This cone, where the norms coincide, must be contained in
$\R_{+}F_{A}$. This implies that $\varphi \in \R_{+}F_{T}$, where $F_{T}$ is some (open)
top--dimensional face of $B_{T}$. The relation $F_{T} \subset F_{A}$ is then a consequence of
McMullen's inequality.

\section{Twisted Alexander polynomials and finite covers} \label{twistal}

\subsection{Twisted Alexander modules and their polynomials} \label{sectionufd}

In this section we are going to define (twisted) Alexander polynomials. These were introduced, for
the case of knots, by Xiao-Song Lin in 1990 (published in \cite{L01}), and his definition was
later generalized to $3$--manifolds by Jiang and Wang \cite{JW93},
Wada
\cite{Wa94}, Kirk-Livingston
\cite{KL99a} and Cha \cite{Ch03}. Our definition below is a generalization
of Cha's approach.

For the remainder of this section let $F$ be a free abelian group and let $R$ be $\Z$ or the
field $\F_p:=\Z/p\Z$ where $p$ is a prime number.  Let $N$ be a compact manifold. Let $\phi \in
Hom(H,F)$ be a non--trivial homomorphism. Through the homomorphism $\phi$, $\pi_{1}(N)$ acts on
$F$ by translations. Let now $A$ be a free finite $R$--module and $\a:\pi_{1}(N) \to \aut_R(A)$
be a representation. We write $A[F]=A\otimes_R R[F]$. We get a representation
\[ \ba{rrcl}  \a\otimes \phi:&\pi_1(N)&\to& \aut(A[F])\\
&g&\mapsto& (\sum_i a_i\otimes f_i \mapsto \sum_i\a(g)(a_i)\otimes (f_i+\phi(g)).\ea \]
We can therefore view $A[F]$ as a left
$\Z[\pi_1(N)]$--module. Note that this module structure commutes with the natural right
$R[F]$--multiplication on $A[F]$.

Let $\ti{N}$ be the universal cover of $N$. Note that $\pi_{1}(N)$ acts on the left on $\ti{N}$
as group of deck transformation. The chain groups $C_*(\ti{N})$ are in a natural way right
$\Z[\pi_1(N)]$--modules, with the right action on $C_{*}(\ti{N})$ defined via $\sigma \cdot g :=
g^{-1}\sigma$, for $\sigma \in C_{*}(\ti{N})$. We can form by tensor product the chain complex
$C_*(\ti{N})\otimes_{\Z[\pi_1(N)]}A[F]$. Now define $H_{*}(N;A[F]):=
H_*(C_*(\ti{N})\otimes_{\Z[\pi_1(N)]}A[F])$, which inherit the structure of right
$\rf$--modules. These modules take the name of twisted Alexander modules.

Our next goal is to define invariants out of these modules. The $\rf$--modules $H_i(N;A[F])$ are
finitely presented and finitely related $\rf$--modules since $\rf$ is Noetherian. Therefore
$H_i(N;A[F])$ has a free $\rf$--resolution
\[ \rf^r \xrightarrow{S_i} \rf^s \to H_i(N;A[F]) \to 0 \]
of finite $\rf$--modules, where we can always assume that $r \geq s$.

\begin{defn} \label{def:alex} The \emph{$i$--th twisted Alexander polynomial} of $(N,\a,\phi)$ is defined
to be  the order of the $\rf$--module $H_i(N;A[F])$,
i.e. the greatest common divisor of the $s\times s$ minors of the $s\times r$--matrix $S_i$. It is denoted by
$\Delta_{N,\phi,i}^{\a}\in \rf$ .
\end{defn}

Note that this definition only makes sense since $\rf$ is a UFD. It is well--known that
$\Delta_{N,\phi,i}^{\a}$ is well--defined only up to multiplication by a unit in $\rf$ and its
definition is independent of the choice of the resolution. Whenever we drop $i$ from
the notation, we will refer to the first Alexander polynomial. If it is clear which manifold $N$ we are referring to,
we will drop it from the notation.

In this paper we will be mostly interested in the case where $F = \Z =\langle t \rangle$ and
$\phi$ is identified with an element of $H^{1}(N)$ (so that $\Delta_{N,\phi,i}\in R[t^{\pm1}]$)
and in the case where $\phi$ is the identity map of $H$. In the latter case we will simply write
$\Delta_{N,i}^{\a}$. Also, we will write ${\Delta}_{N,\phi,i}$  in the case that $\a:\pi_1(N)\to
\gl(\Z,1)$ is the trivial representation. According to these conventions, when $R = \Z$,
$\Delta_{N} \in \Z[H]$ is the ordinary multivariable Alexander polynomial of $N$.

Note that $\Delta_{S^3\sm \nu K}\in \zt$ is just the classical Alexander polynomial $\Delta_K$
of a knot. Furthermore recall that  $H_1(S^{3} \setminus \nu K;\zt)\cong H_1(N(K);\zt)$. In
particular $\Delta_{N(K)}=\Delta_{S^3\sm \nu K}=\Delta_K$.

\subsection{Twisted Alexander polynomials and finite covers}
For the most part of this paper we will be concerned with a particular type of representation:
Let $\a:\pi_1(N) \to G$ be an epimorphism onto a finite group $G$. Then we get an induced
representation $\pi_1(N)\to \aut_R(R[G])$ given by left multiplication. We denote this
representation by $\a$ as well.

To study this case, we first need some definitions and results around the map $\a$. Denote the
normal $G$--cover of $N$ induced by $\a : \pi_{1}(N) \to G$ by $N_G$ (a compact manifold), with
$\pi : N_{G} \rightarrow N$ the covering map. We have an exact sequence \[ 1 \rightarrow
\pi_{1}(N_{G}) \xrightarrow{\pi_*} \pi_{1}(N) \xrightarrow{\a} G \rightarrow 1.\] As $H$ is free
abelian, the inclusion $\pi_{*} : \pi_{1}(N_{G}) \to \pi_{1}(N)$ induces a map $\pi_{*} : H_{G}
\rightarrow H$, where $H_{G}$ is the maximal free abelian quotient of $H_{1}(N_{G})$. We have
therefore a canonical homomorphism $Hom(H,F) \rightarrow Hom(H_{G},F)$ so that given any
homomorphism $\phi: H \rightarrow F$ we can consider the induced homomorphism $\phi_{G} :=
\pi^{*}\phi : H_G \rightarrow F$.

\begin{lemma} \label{lemmamapsurj}
Let $\a:\pi_1(N)\to G$ be an epimorphism to a finite group $G$.
\bn
\item  $\pi^*: Hom(H,F)
\rightarrow Hom(H_{G},F)$ is injective (in particular $b_{1}(N_{G}) \geq b_{1}(N))$.
\item Given $\phi \in Hom(H,F)$ we have
$\mbox{rk}(\im \, \phi) = \mbox{rk}(\im \, \phi_{G})$.
\en
\end{lemma}

\begin{proof}
Recall that we identify $\hom(H,F) = \hom(\pi_{1}(N),F)$ and, similarly, $\hom(H_G,F) =
\hom(\pi_{1}(N_G),F)$. Using our identifications we have
$\pi^{*}\phi = \phi|_{\pi_{1}(N_{G})}$. Let $\phi \in \mbox{ker}\, \pi^{*}$. For any $g \in \pi_{1}(N)$ we have  $\phi(g)^{|G|} =
\phi(g^{|G|}) = 0$ as $g^{|G|} \in \mbox{ker}\, \alpha = \pi_{1}(N_{G})$. But $\phi^{|G|} = 0$
implies $\phi = 0$, as $\hom(\pi_1(N),F)$ is torsion-free. This proves the first statement. The
second statement is immediate.
\end{proof}

Note that $\div \phi_{G} = \div (\phi \circ \pi_{*}) \geq \div \phi \cdot \div \pi_{*}$ and
equality holds by \textit{fiat} if $b_{1}(N) = 1$. We mention also that $\div \pi_{*}$ divides
$|G|$, a result that is useful in interpreting some of the formulae that follow. Proving such
relation is a straightforward exercise in group theory based on Lagrange's formula.

The key tool of this paper is a relation between the twisted Alexander polynomials of $N$
determined by $\a$ and (part of) the ordinary Alexander polynomials for the corresponding
cover $N_{G}$. This relation follows from an isomorphism of Alexander modules that is a
straightforward generalization of the Eckmann-Shapiro Lemma in group cohomology. We have the
following:

\begin{lemma} \label{lemmaalexg}
Let $\phi:H \to F$ be a non--trivial homomorphism to a free abelian group $F$, then
\[ \Delta_{N,\phi,i}^{\a}=\Delta_{N_G,\phi_G,i}\in R[F]. \]
\end{lemma}

\begin{proof}
The proof of this statement follows the steps of the proof of the Eckmann-Shapiro Lemma: first
observe that $\rf$ inherits a left $\Z[\pi_{1}(N_{G})]$--module structure by restriction of its
$\Z[\pi_{1}(N)]$--module structure. Secondly note that the universal covers $\ti{N}$ and
$\ti{N}_G$ of $N$ and $N_{G}$ are the same. Using associativity of the tensor product we have
the following canonical  isomorphisms of chain complexes of $R[F]$--modules
\[ \ba{rcl} C_{*}(\ti{N}_G) \otimes_{\Z[\pi_{1}(N_{G})]}R[F] &=& C_{*}(\ti{N}) \otimes_{\Z[\pi_{1}(N_{G})]} R[F]\\
& \cong& C_{*}(\ti{N}) \otimes_{\Z[\pi_{1}(N)]} (\Z[\pi_{1}(N)] \otimes_{\Z[\pi_{1}(N_{G})]}
R[F])\\
& \cong& C_{*}(\ti{N}) \otimes_{\Z[\pi_{1}(N)]} R[G][F], \ea \] where to establish the last
equality we have used the fact that $\Z[\pi_{1}(N)] \otimes_{\Z[\pi_{1}(N_{G})]} R[F]$ is
canonically isomorphic to $R[G][F]$. This proves the lemma, since $\Delta_{N_G,\phi_G,i}$ is
defined using the first chain complex and $\Delta_{N,\phi,i}^{\a}$ is defined using the last
chain complex.
\end{proof}

As in the ordinary case, the twisted $i$--th Alexander polynomials for $i=0,2,3$ contain only a
modest amount of information.

\begin{lemma} \label{lem:alex02}
Let $N$ be a 3--manifold, $\phi \in Hom(H,F)$ a non--trivial homomorphism to a free abelian group
and $\a:\pi_1(N)\to G$ an epimorphism to a finite group $G$ such that $\Delta_{N,\phi}^\a\ne
0$. Then $\Delta_{N,\phi,3}^\a=1$ and
\[ \Delta_{N,\phi,0}^\a = \Delta_{N,\phi,2}^\a =
\left\{ \begin{array}{ll} 1 & \mbox{if rk$(\im \, \phi_{G}) >
     1$},  \\ \\ (t^{div{\phi_{G}}}-1)
    & \mbox{if $\im \, \phi_{G} = \langle t^{div{\phi_{G}}} \rangle$, $t \in F$ indivisible}. \end{array} \right.
\]
\end{lemma}

\begin{proof} The calculation follows from Lemma \ref{lemmaalexg} and the corresponding well--known
result for the ordinary Alexander polynomial. \end{proof}

\textit{A priori} $\Delta^{\a}_{N,\phi}$ is defined only up to units in $R[F]$. Using Lemma
\ref{lemmaalexg} and a result of Turaev on the first Alexander polynomial of closed 3--manifolds
(\cite[p.~176]{Tu86}) we can actually suppress the ambiguity of the first twisted Alexander
polynomial relative to multiplication by an element of $F$. Therefore we will be able to assume,
whenever suitable, that $\Delta^{\a}_{N,\phi}\in \rf$ is symmetric and well--defined up to a
unit of $R$.

Consider now the complex $C_{*}(N;R[G][F]) = C_*(\ti{N})\otimes_{\Z[\pi_1(N)]}R[G][F]$. Equip
$N$ with a finite CW--structure. Then picking lifts of ordered cells from $N$ to $\ti{N}$ and
using the canonical basis of $R[G]$ as an $R$--module we can view $C_{*}(N;R[G][F])$ as a based
complex of $R[F]$--modules. Let $Q(F)$ denote the field of fractions of $R[F]$. It follows from
Lemma \ref{lem:alex02} and the fact that $Q(F)$ is flat over $R[F]$ that $C_{*}(N;R[G][F])
\otimes_{R[F]} Q(F)$ is acyclic if and only if $\Delta_{N,\phi}^\a\ne 0$.
When this complex is acyclic we can consider its Reidemeister-Franz torsion
$\tau^{\a}_{N,\phi} \in (Q(F)\sm \{0\}) /\{\pm rf|f\in F, r\mbox{ unit in }R\}$. Otherwise we
set $\tau^{\a}_{N,\phi}=0$.  It follows from \cite[Theorem~4.7]{Tu01} together with Lemma
\ref{lem:alex02} that $\tau^{\a}_{N,\phi}\ne 0$ if and only if $\Delta^\a_{N,\phi}\ne 0$,
furthermore if $\Delta^\a_{N,\phi}\ne 0$, then \be \label{equ:taudelta} \tau^{\a}_{N,\phi}   =
\prod_{i=0}^{3} (\Delta^{\a}_{N,\phi,i})^{(-1)^{i+1}} \in Q(F)/\{\pm rf|f\in F, r\mbox{ unit in
}R\}.\ee We denote by $\tau^{\a}_{N}$ the torsion associated to the identity map of $H$, which
corresponds to the multivariable Alexander polynomial (e.g. $\tau_{N}$ is Milnor's torsion). One
of the properties of torsion is its functoriality, encoded in the following theorem:

\begin{theorem} \label{thm:functorial} Let $\rho : R \to \ti{R}$ be a ring homomorphism and
let $\phi$ be a non--trivial element of $Hom(H,F)$. These induce a ring
homomorphism $\rho \otimes \phi :R[H]\to \ti{R}[F]$. Now let \[ Q=\{ fg^{-1} | f,g \in R[H],
(\rho\otimes \phi)(g)\ne 0\} \subset Q(R[H]).\] Then $\tau^{\a}_{N,\phi}\in Q$. Furthermore $\rho\otimes \phi$ extends to a homomorphism $Q\to Q(\ti{R}[F])$ and
\[ (\rho \otimes \phi)(\tau^{\a}_{N}) = \tau^{\a}_{N,\phi}, \]
this equality taking value in $\ti{R}[F]$ for $b_{1}(N) \geq 2$ and in ${Q}(\ti{R}[F])$ for
$b_{1}(N) = 1$.
\end{theorem}

\begin{proof} The proof follows again from Lemma \ref{lemmaalexg} and functoriality for the ordinary torsion, proved
in \cite[Theorem~1.1.3]{Tu86}. (Note that such result is stated therein only for $R =
\ti{R} = \Z$, but the proof carries over to the general case.)
\end{proof}

\noindent The following is an immediate corollary to Lemma \ref{lem:alex02} and Theorem
\ref{thm:functorial}.

\begin{proposition} \label{relalex} Let $N$ be a  $3$--manifold and let
$\a: \pi_{1}(N) \rightarrow G$ be an epimorphism onto a finite group. Let $\pi_{*} : H_G \to H$ be
the induced map. Then the twisted Alexander polynomials of $N$ and the ordinary Alexander
polynomial of $N_{G}$ satisfy the following relations: \\ If $b_{1}(N_{G}) > 1$, then
\begin{equation}  \label{proone} \Delta^{\a}_{N}   = \left\{ \begin{array}{ll} \pi_{*}(\Delta_{N_{G}})  & \mbox{if $b_{1}(N) >
     1$}, \\ \\ (t^{div{\pi_{*}}}-1)^{2}\pi_{*}(\Delta_{N_{G}}) & \mbox{if $b_{1}(N) = 1$,
$\im \, \pi_{*} = \langle t^{div\, \pi_{*}} \rangle$, $t \in H$ indivisible}.
\end{array} \right.
\end{equation} If $b_{1}(N_{G}) = 1$, then $b_{1}(N) = 1$ and
\begin{equation} \label{protwo} \Delta^{\a}_{N} = \pi_{*} (\Delta_{N_{G}}). \end{equation} \end{proposition}

\begin{proof} First assume that $\pi_*(\Delta_{N_G})=0$. Distinguishing the cases where $\Delta_{N_G} = 0$ and
$\Delta_{N_G} \neq 0$, and applying in the latter case Lemma \ref{lem:alex02} and Equation
(\ref{equ:taudelta}), this implies that $\pi_*(\tau_{N_G})=0$.
It follows from Theorem \ref{thm:functorial}, applied to $\pi_{*} \in
Hom(H_{G},H)$, that $\tau_{N_G,\pi_*}=0$, hence $\Delta_{N_G,\pi_*}=\Delta_N^\a=0$. Now assume
that $\pi_*(\Delta_{N_G})\ne 0$; then we have by Theorem \ref{thm:functorial} and Equation
(\ref{equ:taudelta}) that
\[ \pi_{*} \Big( \prod_{i=0}^{3} (\Delta_{N_{G},i})^{(-1)^{i+1}} \Big)
= \prod_{i=0}^{3} (\Delta_{N_{G},\pi_{*},i})^{(-1)^{i+1}}. \]
Apply now Lemma \ref{lemmaalexg} choosing as $\phi$ the identity map of $H$ (observing that, for
this choice, $\phi_G=\pi_*$) and apply Lemma \ref{lem:alex02} twice. Using Lemma
\ref{lemmamapsurj}, which implies that $b_1(N)=\mbox{rk}(\im \, \pi_{*})$, Equations
(\ref{proone}) and (\ref{protwo}) follow.
\end{proof}

The proposition above asserts that the twisted Alexander polynomial of $N$ is determined by the
Alexander polynomial of $N_{G}$ or, stated otherwise, by the Seiberg-Witten invariants of
$N_{G}$.

\subsection{Twisted Alexander norms and Thurston norm} \label{section:twinorm}

We can introduce a norm on $H^{1}(N;\R)$ determined by the (symmetric) twisted Alexander
polynomial $\Delta_{N}^{\a}$, by mimicking McMullen's definition of the Alexander norm.
Precisely we have the following:
\begin{defn} Let $\Delta_{N}^{\a}$ be the twisted Alexander polynomial associated to an epimorphism
$\a: \pi_{1}(N) \rightarrow G$. The \emph{twisted Alexander norm} associated to $\a$ is the
(semi)norm on $H^{1}(N;\R)$ defined as \[ \|\phi\|_{A}^{\a}:= 2 \, \mbox{max} \{ \phi(g) | g \in
\mbox{supp} \Delta_{N}^{\a} \subset H \}. \] As for the ordinary case, this is understood to
mean that $\|\phi\|_{A}^{\a} = 0$ if $\Delta_{N}^{\a} = 0$.
\end{defn} This norm satisfies, from its definition, all the usual properties of the ordinary
Alexander norm. In particular, its unit ball $B_{A}^{\a}$ is a convex (possibly non--compact)
polyhedron, dual to the Newton polyhedron of $\Delta_{N}^{\a}$. Whenever useful, we will specify
in the notation the manifold to which the norm refers to.

Let $\phi \in H^{1}(N)$. We can think of $\phi$ as an element of the integer lattice of
$H^{1}(N;\R)$. As in the classical Alexander theory, there is a relation between the degree of
$\Delta^{\a}_{N,\phi}$ and the twisted Alexander norm of $\phi$. (Note that for $f=\sum_{i=n}^m
a_it^i, a_na_m\ne 0$ we define $\deg(f)=m-n$, furthermore we write $\deg(0)=-\infty$ with the
understanding that $-\infty < a$ for any $a\in \Z$.)
The latter, moreover, is related
to the ordinary Alexander norm of $\phi_{G}$ on $N_{G}$. The form of these relations is
expressed in the following proposition.
\begin{proposition} \label{holdt}
Let $\phi$ be a non--trivial element of $H^{1}(N)$; then the following inequalities hold true:
\begin{equation} \label{threelines} \deg\Delta^{\a}_{N,\phi}   \leq \| \phi \|_{A,N}^{\a}  +
\left\{ \begin{array}{l} 2 \div \, \phi_{G} \\ \\ 0 \\ \\ 0 \end{array} \right. \leq \| \phi_{G}
\|_{A,N_{G}} + \left\{ \begin{array}{ll} 2 \div \, \phi_{G} & \mbox{if $b_{1}(N_{G}) \geq
b_{1}(N)
>
     1$}, \\ \\ 2 \div \, \phi_{G} & \mbox{if $b_{1}(N_{G}) > b_{1}(N) =
     1$}, \\ \\ 0 & \mbox{if $b_{1}(N_{G}) = b_{1}(N) =
     1$}. \end{array} \right. \end{equation}
The first relation (provided that $\Delta^{\a}_{N} \neq 0$) is an equality for all $\phi$ in the
cone over the top--dimensional faces  of the unit ball of the twisted Alexander norm (in
particular, whenever $b_{1}(N) = 1$).
\end{proposition}

\begin{proof} Theorem \ref{thm:functorial}, applied to $\phi \in \hom(H,\Z)$ asserts that, denoting
$\Delta_{N,\phi}^{\a} = \sum a_{i} t^{i}$ and $\Delta_{N}^{\a} = \sum_{g \in H} b_{g} g$, we
have
\begin{equation} \label{assert} \sum a_{i} t^{i} = \left\{ \begin{array}{ll}
(t^{div{\phi_{G}}}-1)^{2}\phi(\Delta^{\a}_{N}) = (t^{div \phi_{G}} - 1)^{2} \sum_{g \in H} b_{g}
t^{\phi(g)}  & \mbox{if $b_{1}(N) >
     1$}, \\ \\ \phi(\Delta^{\a}_{N}) = \sum_{g \in H} b_{g} t^{\phi(g)}  & \mbox{if $b_{1}(N) = 1$.} \end{array} \right.
\end{equation}
This implies the first inequality. Whenever $\phi$ is in the cone over a top--dimensional face
of the unit ball of the twisted Alexander norm it attains its maximum value at a single $g\in
\mbox{supp} \Delta_N^\a\subset H$ which is dual to that face. In that case, the equality holds.
For the second inequality, we can use Proposition \ref{relalex} to determine the relationship
between the twisted Alexander norm of $N$ and the Alexander norm of $N_{G}$. Write
$\Delta_{N_{G}} = \sum_{g \in H_{G}} c_{g} g \in R[H_{G}]$. Let us consider first the cases
where $b_{1}(N_{G}) \geq b_{1}(N) > 1$ or $b_{1}(N_{G}) = b_{1}(N) = 1$: using Proposition
\ref{relalex}, $\Delta_{N}^{\a} = \sum_{g \in H_G} c_{g} (\pi_{*} g) \in R[H]$; it follows that
the support of $\Delta^{\a}_{N}$ is contained in the image of the support of $\Delta_{N_{G}}$,
and when $b_{1}(N_{G}) > b_{1}(N)$ this inclusion may happen to be strict. We deduce that, for
all $\phi \in H^{1}(N;\R)$,
\[ \|\phi \|_{A,N}^{\a} = 2 \, \mbox{max} \{\phi(\pi_{*}g) | \pi_{*} g \in \mbox{supp} \Delta_{N}^{\a} \}
\leq  2 \, \mbox{max} \{ \phi_{G}(g) | g \in \mbox{supp} \Delta_{N_{G}} \} =
\|\phi_{G}\|_{A,N_{G}}. \] This gives the first and third line of Equation (\ref{threelines}).
In the remaining case where $b_{1}(N_{G}) > b_{1}(N) = 1$ Proposition \ref{relalex} gives
$\Delta_{N}^{\a} = (t^{div \pi_{*}} - 1)^{2} \sum_{g \in H_G} c_{g} (\pi_{*} g) \in R[H]$ which
implies that, for all $\phi \in H^{1}(N)$,
\[ \ba{rcl} \|\phi\|_{A,N}^{\a} &= &2 \div \pi_{*} \cdot \div \phi + 2 \mbox{max} \{ \phi(\pi_{*}g) | \pi_{*} g \in
\mbox{supp} \Delta^{\a}_{N}  \}   \\[0.2cm] &\leq& 2 \div \phi_{G} +  2 \mbox{max} \{ \phi_{G}(g) |
g \in \mbox{supp} \Delta_{N_{G}} \}\\[0.2cm] &=& 2 \div \phi_{G} + \|\phi_{G}\|_{A,N_{G}}.
\ea \] This completes the second line of the inequality and finishes the proof.
\end{proof}

We can now examine the relationship between the twisted Alexander norm of $N$ and the Thurston
norm, to obtain an inequality \`a la McMullen. We have the following Lemma:

\begin{lemma} \label{normineq} Let $\phi \in H^1(N)$ be non--trivial. Then we have the following inequality:
\begin{equation} \| \phi \|^{\a}_{A}  \leq |G| \cdot \| \phi \|_{T} + \left\{ \begin{array}{ll} 0 & \mbox{if $b_{1}(N) >
     1$},  \\ \\ 2 \div \phi_{G}  & \mbox{if $b_{1}(N) =
     1$}. \end{array} \right. \end{equation} \end{lemma}

\begin{proof} Corollary 6.13 of \cite{Ga83} states that
\[ |G|\cdot||\phi||_{T,N}= ||\phi_G||_{T,N_G}.\] Combined with McMullen's inequalities for $N_{G}$ and
(\ref{threelines}) this gives the desired inequalities. \end{proof}

\section{Main results} \label{fini}

%============================

\subsection{Finite covers and fibered manifolds}

Let $N$ be a 3--manifold and let $\a:\pi_1(N)\to G$ be an epimorphism to a finite group $G$ and
denote as above by $\pi: N_{G} \rightarrow N$ the induced covering map. Assume that $(N,\phi)$
fibers over $S^{1}$ for $\phi \in H^{1}(N)$: then it is well known that $(N_G,\phi_{G})$ fibers
over $S^{1}$ too. (In terms of differential forms, the fact that $\phi$ is fibered corresponds
to the fact that there exists a non--degenerate closed $1$--form $\nu \in \Omega^{1}(N,\R)$,
unique up to isotopy, representing the corresponding class in $H^{1}_{DR}(N;\R)$. The
representative for $\phi_{G}$ is then given by $\pi^{*} \nu$.) The Euler class of the fibration
$\phi_{G}$, as element of $H_{G}$, is given by $e_{F_{T,N_{G}}} = PD_{N_{G}} (\pi^{*} PD_{N}
(e_{F_{T,N}}))$.

We want to analyze the properties of the twisted Alexander polynomial of a fibered
$3$--manifold. For the rest of the section, we are going to fix $R = \Z$.

We can apply Theorem \ref{thmmcmullen} and Proposition \ref{relalex} to $N_{G}$ to get a
constraint on the twisted Alexander and Thurston norms of a fibered class, that is a ``twisted
version" of Theorem \ref{thmmcmullen}.  Precisely, we have the following:

\begin{theorem} \label{mainthmg} Let $N \neq S^{1} \times S^{2}$ be a fibered $3$--manifold,
$F_{T}$ be a fibered face of the Thurston unit ball, and let $\a:\pi_1(N)\to G$ be an
epimorphism onto a finite group. Then the following hold true:
\bn
\item For all $\phi \in \R_{+}F_{T} \cap H^{1}(N)$ \[ \| \phi \|^{\a}_{A}  = 2 \phi(\lambda_{G}) =
|G| \| \phi \|_{T} + \left\{ \begin{array}{ll} 0 & \mbox{if $b_{1}(N) >
     1$},  \\ \\ 2 \div \, \phi_{G}  & \mbox{if $b_{1}(N) =
     1$}, \end{array} \right. \] where $\lambda_{G}$ is a vertex of the Newton polyhedron of $\Delta^{\a}_{N}$ with
     coefficient $\pm 1$.
If $b_{1}(N) > 1$, then $\lambda_{G} = -\frac{|G|}{2}e_{F_{T}} = |G|\lambda$ and $F_{T} \subset
|G|F^{\a}_{A}$ (where $F^{\a}_{A}$ is the top--dimensional face of the Alexander unit ball dual to
$\lambda_{G}$); if $b_{1}(N) = 1$, then $\lambda_{G} = -\frac{|G|}{2}e_{F_{T}} + \div \, \pi_{*} =
|G|\lambda + \div \, \pi_{*} - |G|$.
\item For all $\phi \in \R_{+}F_{T} \cap H^{1}(N)$,
$\Delta_{N,\phi}^{\a} \in \Z[t^{\pm 1}]$ is monic and $\deg\, \Delta_{N,\phi}^{\a} = |G|
\, \|\phi\|_{T} + 2 \div \, \phi_{G}$.
\en
\end{theorem}

\begin{proof} Write $\Delta_{N_G}=\sum_{g \in H_{G}} c_gg$. As $(N_{G},\phi_{G})$ fibers,
Theorem \ref{thmmcmullen} states that the Newton polyhedron of $\Delta_{N_{G}}$
has a vertex $\lambda$ such that $c_{\lambda} = \pm 1$, dual to $\phi_{G}$. Define
\begin{equation} \label{lambda} \lambda_{G} = \pi_{*} \lambda +  \left\{ \begin{array}{ll} 0 & \mbox{if $b_{1}(N_{G}) \geq b_{1}(N) >
     1$ or $b_{1}(N_{G}) = b_{1}(N) =
     1$}, \\ \\ \div \, \pi_{*} & \mbox{if $b_{1}(N_{G}) > b_{1}(N) =
     1$}.  \end{array} \right. \end{equation}
We claim that $\lambda_{G}$ is a vertex of the Newton polyhedron of $\Delta_{N}^{\a}=\sum_{g \in
H} b_gg$ such that $b_{\lambda_{G}} = \pm 1$. In fact by Theorem \ref{thmmcmullen} $\phi(\pi_{*}
\lambda) = \phi_{G}(\lambda)$ is the \textit{unique} maximum over all elements in $\mbox{supp}
\Delta_{N_{G}} \subset H_{G}$, which entails that $\lambda_{G}$ is a vertex of $\mbox{supp}
\Delta^{\a}_{N}$, dual to $\phi$. Proposition \ref{relalex} then implies that $b_{\lambda_{G}} =
c_{\lambda} = \pm 1$.

As stated in Theorem \ref{thmmcmullen} the vertex $\lambda$ equals $-\frac12 e_{F_{T,N_{G}}}$ if
$b_{1}(N_{G}) > 1$ and $-\frac12 e_{F_{T,N_{G}}} + 1$ if $b_{1}(N_{G}) = 1$. A straightforward
calculation shows that $\pi_{*} e_{F_{T,N_{G}}} = |G| e_{F_{T,N}}$, which in turn implies the
relation between $\lambda_{G}, e_{F_{T}}$ and $\lambda$ stated in (1). This, together with the
equality $\varphi(-e_{F_{T}}) = \|\varphi\|_{T}$ that holds for all $\varphi \in \R_{+}F_{T}$,
implies the equalities of (1).

The second part follows either from  Theorem \ref{thmmcmullen} (2) and Proposition \ref{relalex} or
it can be deduced from (1) as follows. Applying Theorem \ref{thm:functorial} as in Equation
(\ref{assert}), we can write
\begin{equation}  \Delta^{\a}_{N,\phi}   = \left\{ \begin{array}{ll} (t^{div{\phi_{G}}}-1)^{2}\phi(\Delta^{\a}_{N}) =
(t^{div{\phi_{G}}}-1)^{2}(\sum_{g \in H} b_{g} t^{\phi(g)}) & \mbox{if $b_{1}(N) > 1$,} \\ \\
\phi(\Delta^{\a}_{N}) = \sum_{g \in H} b_{g} t^{\phi(g)} & \mbox{if $b_{1}(N) = 1$}. \end{array}
\right. \end{equation} In either case, by (1), the condition that $\phi$ is contained in
$\R_{+}F_{T} \subset \R_{+}F_{A}^\a$ entails that $\phi(g)$ attains maximal value in correspondence
of the unique $\lambda_{G} \in \mbox{supp} \Delta^{\a}_{N}$, with coefficient $b_{\lambda_{G}} =
\pm 1$. It follows that $\Delta^{\a}_{N,\phi}$ is monic and
\begin{equation}
\deg\Delta^{\a}_{N,\phi}   = \left\{ \begin{array}{ll} 2 \phi (\lambda_{G}) +
2 \div \phi_{G} = |G|\|\phi\|_{T} + 2 \div \phi_{G} & \mbox{if $b_{1}(N) > 1$}, \\ \\
2 \phi (\lambda_{G}) = |G|\|\phi\|_{T} + 2 \div \phi_{G} & \mbox{if $b_{1}(N) = 1,$}
\end{array} \right. \end{equation}
where for the last equality we have used part (1) and the relation $\div \phi_{G} = \div \pi_{*}
\cdot \div \phi$.
\end{proof}

It is important to point out that the constraint of part (2) of Theorem \ref{mainthmg} (and,
later, its symplectic counterpart) is, in practice, the condition that can be directly tested to
investigate if a class $\phi$ represents a fibration.

Our results in Section \ref{sectionapp} and Section \ref{section:ex} (cf. also
\cite{FK06,Ch03,GKM05}) show that twisted Alexander polynomials are very efficient at
detecting non--fibered manifolds.  This leads us to conjecture the following:

\begin{conjecture} \label{conjcha}
Let $N$ be a $3$--manifold and let $\phi \in H^{1}(N)$ a non--trivial class such that for
\textit{any} epimorphism onto a finite group $\a: \pi_{1}(N) \rightarrow G$ the twisted
Alexander polynomial $\Delta^{\a}_{N,\phi} \in \Z[t^{\pm 1}]$ is monic and
$\deg\Delta_{N,\phi}^{\a} = |G| \, \|\phi\|_{T} + 2 \div \, \phi_{G}$. Then $(N,\phi)$
fibers over $S^{1}$. \end{conjecture}

Note that this conjecture is slightly stronger than the converse of Theorem \ref{mainthmg}, as
we are not requiring \textit{a priori} that $\phi$ is in the cone over a top-dimensional face of
the Thurston unit ball. In the next subsection, we will discuss how this conjecture relates to
Conjecture \ref{conjfolk}. Further evidence to this conjecture is discussed in Section
\ref{sectionconj}.\\

\emph{Added in proof.} In \cite{FV07} Conjecture \ref{conjcha} has been proved in the two special cases that $N$ has vanishing Thurston norm or that $N$ is a graph manifold.

%===========================================
%===========================================
\subsection{Finite covers and symplectic manifolds}
We can use the properties of Theorem \ref{mainthmg} to further analyze the relation between
fibered 3--manifolds and symplectic $S^{1} \times N$, to get more evidence to
Conjecture \ref{conjfolk}. The proof of Theorem \ref{mainthmg} boils down to the fact that finite covers of
fibered manifolds are fibered again. Similarly, we will use the fact that finite covers of a
symplectic manifold admit a symplectic structure to obtain our main theorem.

\begin{theorem} \label{strongsymp} \label{mainthm}
Let $N$ be an irreducible $3$--manifold such that $S^{1} \times N$ admits a symplectic structure
with K\"unneth component $\varphi \in H^{1}(N)$ and canonical class $K$. Then there exists a
top--dimensional face $F_{T}$ of the Thurston unit ball such that $\varphi \in \R_{+}F_{T}$ for
which the same properties as in (1) and (2) of Theorem \ref{mainthmg} hold true, with
$\lambda_{G}$ equal to $\frac{|G|}{2} PD_{N} K$ if $b_{1}(N) > 1$ and $\frac{|G|}{2} PD_{N} K +
\div \, \pi_{*}$ if $b_{1}(N) = 1$. \end{theorem}

\begin{proof}
First, observe that, denoting (with slight abuse of notation) the finite cover $S^{1} \times
N_{G} \rightarrow S^{1} \times N$ by $\pi$, the manifold $S^{1} \times N_{G}$ is symplectic with
symplectic form $\pi^{*} \w$ whose cohomology class has K\"unneth component $\pi^{*} \varphi =
\varphi_{G}$, and canonical class $\pi^{*}K$ (that we can consider as an element of
$H^{2}(N_{G})$). Moreover, the assumption of irreducibility of $N$ implies that $N_{G}$ is
irreducible too, by the Equivariant Sphere Theorem (cf. \cite{Du85}). Therefore we can apply
Theorem \ref{thmvidussi} for  $N_{G}$ in the same way we applied Theorem \ref{thmmcmullen} in
the proof of Theorem \ref{mainthmg}, with $\lambda$ equal to $\frac12 PD_{N_{G}} \pi^{*} K$ if
$b_{1}(N_{G}) > 1$ and equal to $\frac12 PD_{N_{G}} \pi^{*} K + 1$ if $b_{1}(N_{G}) = 1$. The
class $\lambda_{G}$, defined as is Equation (\ref{lambda}), is a vertex of the Newton polyhedron
of $\Delta_{N}^{\a}$, dual to $\varphi$, with coefficient $b_{\lambda_{G}} = \pm 1$. The relation
between the Thurston norm and the twisted Alexander norm of $\varphi$ follows {\it verbatim}
that of Theorem \ref{mainthmg}, with $K$ playing the role of $-e_{F_{T,N}}$. To extend the
result to the entire cone $\R_{+}F_{T}$, when $b_{1}(N) \geq 2$, we use the inclusion of the
unit balls $B_{T} \subset |G|B_{A}^{\alpha}$ to grant that, as  $F_{T}$ and $|G|F_{A}^{\a}$
intersect, $F_{T} \subset |G|F_{A}^{\a}$. Finally, to obtain (2), we proceed exactly as in the
proof of Theorem \ref{mainthmg}.
\end{proof}

For sake of clarity, we give a slightly weaker form of (2) of Theorem \ref{strongsymp} that will
provide, in what follows, the obstruction to the existence of symplectic forms which we mostly
use.

\begin{proposition} \label{sympcorg}
Let $N$ be an irreducible $3$--manifold such that $S^{1} \times N$ admits a symplectic structure
with K\"unneth component $\phi \in H^{1}(N)$.  Let $\a:\pi_1(N)\to G$ be an epimorphism onto a
finite group, then $\Delta_{N,\phi}^\a \in \Z[t^{\pm 1}]$ is monic and
\[ \deg\, \Delta_{N,\phi}^\a = |G| \, \|\phi\|_{T} + 2 \div \, \phi_{G}.\]
\end{proposition}

\begin{remark} \bn
\item  Note that (2) of Theorem \ref{mainthmg} (respectively Proposition \ref{sympcorg}) implies that in order
for a class $\phi$ to be fibered (respectively for $S^{1} \times N$ to admit a symplectic form
with K\"unneth component $\phi$) the degrees of the twisted Alexander polynomials must satisfy
$\deg \Delta_{N,\phi}^{\a} =|G| \deg \Delta_{N,\phi} + 2(\div \pi_{*} - |G|) \div \phi$. This
condition can be verified, for a given choice of $\a$, regardless of our explicit knowledge of
the Thurston norm of $\phi$, on which we may have no information. In this sense, these results
have {\it a priori} a larger field of application than Theorem \ref{thmmcmullen}.
\item
The content of Theorem \ref{strongsymp} follows, ultimately, from the symplectic constraints on
Seiberg-Witten invariants of finite covers of $S^{1} \times N$. In this sense, this result is
just an application of the well--known principle that we can unveil information about a manifold
by looking at the Seiberg-Witten invariants of its finite covers.
\item It is not difficult to verify that, when $G = \Z_{d}$ and $\a: \pi_{1}(N) \to H \to \Z_{d}$ (i.e. when $N_{G}$ is a cyclic cover), Theorem \ref{strongsymp} does not provide any new obstruction, as in that case the twisted Alexander polynomial is completely determined in terms of the ordinary Alexander polynomial.
\item In \cite{Ba01}, the author asks whether the orbit space of any symplectic $4$-manifold with a free $S^{1}$--action is fibered; an affirmative answer would include, in particular, Conjecture \ref{conjfolk}.  It is possible that proceeding along the lines of Theorem \ref{strongsymp} a similar statement can be obtained.
\en
\end{remark}

In the next sections we will use these results to show the presence of new obstructions to the
existence of symplectic structure on certain $S^{1} \times N$'s for which the conditions of
Theorem \ref{thmvidussi} are instead satisfied. But, perhaps more importantly, we can observe
that analyzing the statement of Theorem \ref{strongsymp} (or Proposition \ref{sympcorg}) we see
that an irreducible $N$ for which $S^{1} \times N$ is symplectic satisfies the conditions of
Conjecture \ref{conjcha}. We have therefore the following:
\begin{corollary} Conjecture
\ref{conjcha} implies Conjecture \ref{conjfolk} for irreducible $3$--manifolds.
\end{corollary}

The condition of irreducibility can actually be removed assuming that the Geometrization
Conjecture holds true, see Proposition \ref{lem:prime}. One appeal of Conjecture \ref{conjcha}
is that it reduces Conjecture \ref{conjfolk} to a purely $3$--dimensional question (see also Section
\ref{sectionconj}), that, although possibly stronger, has a distinctively more ``quantitative"
flavor than Conjecture \ref{conjfolk}. Note that, in fact, Conjecture \ref{conjfolk} follows also from the slightly
weaker condition that Conjecture \ref{conjcha} holds for the classes contained in a cone over
the top--dimensional faces of the Thurston unit ball. In both cases Conjecture \ref{conjfolk}
is implied in the strongest sense, namely Conjecture \ref{conjcha} implies that the
restriction of the symplectic cone(s) of $S^{1} \times N$ to $H^{1}(N;\R)$ coincides with the
fibered cone(s) of $N$.

\section{Abelian invariants do not determine the existence of symplectic structures}
\label{sectionapp}

In this section we will exhibit a family of examples where, applying Theorem \ref{strongsymp},
we will be able to exclude the existence of symplectic structures, even though the ordinary
Alexander polynomial does not provide any obstructions. These examples, in a sense, are defined
\textit{ad hoc} to allow a simple computation of the twisted Alexander polynomial for suitable
$G$. More challenging, and perhaps more interesting examples are discussed in Section
\ref{section:ex}.

The examples we consider are obtained by $0$-surgery on an infinite family of non--fibered knots
$K_{i}$. Application of Theorem \ref{strongsymp} is made possible by the fact that such
manifolds, by \cite{Ga87b}, are irreducible. Actually, these examples allow us to show a
stronger result, namely that abelian invariants of $K$ can not determine the existence of
symplectic structures on a manifold of the form $S^1\times N(K)$. We have the following:

\begin{theorem} \label{thmsamev} Let $K$ be a non--trivial fibered knot. Then there exist infinitely many
distinct knots $K_i,i\in \N$ such that for any $i$
\bn
\item $\Delta_{K_i}=\Delta_K$, in particular $\Delta_{K_i}$ is monic,
\item $\deg(\Delta_{K_i})=2\, g(K_i)$,
\item $S^1\times N(K_i)$ is not symplectic,
\item furthermore $N(K_i)\ne N(K_j)$ for any $i\ne j$.
\en
\end{theorem}

The proof is somewhat similar to an argument of Cha in \cite{Ch03}
concerning the existence of non--fibered knots with certain properties, and it will require the
remainder of this section. Before proving this theorem, we briefly summarize the satellite
construction for knots.

Let $K,C$ be knots. Let $A\subset S^3\sm \nu K$ be a curve, unknotted in $S^3$. Then $S^3\sm \nu
A$ is a solid torus. Let $h:\partial(\overline{\nu A})\to \partial(\overline{\nu C})$ be a
diffeomorphism which sends a meridian of $A$ to a longitude of $C$, and a longitude of $A$ to a
meridian of $C$. The space
\[ \left({S^3\sm \nu A}\right) \cup_{h} \left({S^3\sm \nu C}\right)\]
is diffeomorphic to $S^3$. The image of $K$ is
denoted by $S=S(K,C,A)$ and called a satellite knot.
Note that we replaced a tubular neighborhood of $C$ by a
knot in a solid torus, namely by $K\subset {S^3\sm \nu A}$.

It is easy to see that $N(S)=\left({N(K)\sm \nu A}\right) \cup_h \left({S^3\sm \nu C}\right)$.
By mapping $S^3\sm \nu C$ to $\overline{\nu A}$ we get an induced map $N(S)\to N(K)$.
Now let $\a:\pi_1(N(K))\to G$ be a homomorphism to a finite group $G$; by the above we get
an induced map $\pi_1(N(S))\to G$. Note that $A$ defines an element in $[A]\in \pi_1(N(K))$
which is well--defined up to conjugation. Let $r$ be the order of $\a([A])\in G$.  Denote by
$L(C)_r$ the $r$--fold cyclic cover of $S^3$ branched along $C$.

\begin{lemma} \label{lemmahofsat} Assume that $[A]=0\in H_1(S^{3} \setminus \nu K)$. Then there exists a short exact
sequence \[ 0\to H_1(L(C)_r)^{\frac{|G|}{r}}\otimes_{\Z} \zt \to H_1(N(S);\zgt)\to
H_1(N(K);\zgt)\to 0.\]
\end{lemma}

\begin{proof}
Consider the following commutative diagram of Mayer--Vietoris sequences of homology groups
(where $R=\zgt$):

\[ \ba{cccccccccc}
\to \hspace{-0.2cm}& H_1(\partial(\overline{\nu A});R) & \hspace{-0.2cm}\to
\hspace{-0.2cm}&H_1({N(K)\sm \nu A};R) &\hspace{-0.2cm}\oplus \hspace{-0.2cm} & H_1({S^3\sm \nu
C};R)& \hspace{-0.2cm}\to \hspace{-0.2cm}&
H_1(N(S);R)&\hspace{-0.2cm}\to\hspace{-0.2cm}& 0 \\[0.1cm]
&\downarrow \cong &&\downarrow\cong &&\downarrow f&&\downarrow&\\[0.1cm]
\to \hspace{-0.2cm}& H_1(\partial(\overline{\nu A});R) &\hspace{-0.2cm}\to\hspace{-0.2cm}
&H_1({N(K)\sm \nu A};R) &\hspace{-0.2cm}\oplus\hspace{-0.2cm}& H_1(\overline{\nu A};R)
&\hspace{-0.2cm}\to \hspace{-0.2cm}& H_1(N(K);R)&\hspace{-0.2cm}\to\hspace{-0.2cm}& 0.\ea
\]
Note that $[A]=0\in H_{1}(S^{3} \setminus \nu K)$ implies that $H_1(S^3\sm \nu C)\to H_1(N(S))$
is the zero map; it follows that
\[ \ba{rcl} H_1(S^3\sm \nu C;\zgt)&\cong
&H_1((S^{3}
\setminus \nu C)_r)^{\frac{|G|}{r}}\otimes_{\Z} \zt, \\
H_1(\overline{\nu A};\zgt)&\cong&\Z^{\frac{|G|}{r}}\otimes_{\Z} \zt,\ea \] where $(S^{3}
\setminus \nu C)_r$ denotes the $r$--fold cyclic cover of $S^{3} \setminus \nu C$. Furthermore
$H_1((S^{3} \setminus \nu C)_r)=H_1(L(C)_r)\oplus \Z$,
where the $\Z$--summand is generated by the lift of $r$--copies of the meridian of $C$. It
follows that \[ \ker\left\{H_1({S^3\sm \nu C};\zgt)\to H_1(\overline{\nu A};\zgt)\right \} \,
\cong \, H_1(L(C)_r)^{\frac{|G|}{r}}\otimes_{\Z} \zt.\] The proof now follows easily from the
above commutative diagram.
\end{proof}

\noindent We are now in a position to prove Theorem \ref{thmsamev}.

\begin{proof}[Proof of Theorem \ref{thmsamev}] Let $K$ be a non--trivial fibered knot with fiber $F$. Let $A\subset S^{3} \setminus \nu K$ be a
simple closed curve disjoint from $F$ that, after surgery, represents a non--trivial element in
$\pi_1(N(K))^{(1)}/\pi_1(N(K))^{(2)}$. (Recall that the derived series of a group $\pi$ is
defined inductively via $G^{(0)}=G, G^{(i+1)}=[G^{(i)},G^{(i)}]$.) Note that we can find such an
$A$ since $\pi_1(S^3\sm K)\to \pi_1(N(K))$ is surjective and since
$\pi_1(N(K))^{(1)}/\pi_1(N(K))^{(2)}\cong H_1(N(K);\zt)\cong \Z^{2genus(K)}\ne 0$. As the
linking number of $A$ and $K$ is zero, $[A]$ is trivial in $H_{1}(N(K))$. After local crossing
changes if necessary we can assume that $A$ is unknotted in $S^{3}$. As $N(K)$ is irreducible
and since $b_{1}(N(K)) =1$ it follows that $N(K)$ is Haken, hence $\pi_1(N(K))$ is residually
finite (cf. \cite{He87}). We can therefore find a finite group $G$ and an epimorphism
$\a:\pi_1(N(K))\to G$ such that $\a([A])\ne e$. Denote by $r$ the order of $\a(A)$ in $G$. Note
that $r>1$.

Let $C_i$ to be connected sum of $i$ copies of the figure--8 knot and of $i$ copies of the
trefoil knot. It follows from \cite[Lemma~4.4]{Ch03} that $H_1(L(C_{i})_r)$ is non--trivial.

From Lemma \ref{lemmahofsat} we have $H_1(N(K_i);\zt)\cong H_1(N(K);\zt)$, in particular
$\Delta_{K_i}=\Delta_K$. Also note that the image of $F$ in $\left({S^3\sm \nu A}\right)
\cup_{\phi} \left({S^3\sm \nu C_i}\right)$ is still a Seifert surface for $K_i$ since $A\subset
S^3\sm F$, hence from Theorem \ref{thmmcmullen} it follows that $2
g(K_i)=2g(K)=\deg(\Delta_{K_i})$.

We now show that $S^1\times N(K_i)$ is not symplectic. Let $\phi\in H^1(N(K_i))$ be one of the two
generators. By Lemma \ref{lemmahofsat} there exists a short exact sequence
\[ 0\to H_1(L(C_i)_r)^{\frac{|G|}{r}}\otimes_{\Z} \zt\to H_*(N(K_i);\zgt)\to H_1(N(K);\zgt)\to 0.\]
Given a finitely generated $\zt$--module $H$ we can define its order $\ord(H)$ as in Definition
\ref{def:alex}.  By the multiplicativity of the order of $\zt$--modules (cf.
\cite[Lemma~5,~p.~76]{Le67}) we have
\[ \Delta_{N(K_i)}^\a=\Delta_{N(K)}^\a\, \ord\left(H_1(L(C_i)_r)^{\frac{|G|}{r}}\otimes_{\Z} \zt\right).\]
It is easy to see that
\[\ord\left(H_1(L(C_i)_r)^{\frac{|G|}{r}}\otimes_{\Z}
\zt\right)=|H_1(L(C_i)_r)|^{\frac{|G|}{r}}.\] It follows that $\Delta_{N(K_i)}^\a$ is not monic,
hence $S^1\times N(K_i)$ is not symplectic by Theorem \ref{mainthm}.

It remains to show that $N(K_i)\ne N(K_j)$ for $i\ne j$. Various methods are viable. We will
distinguish the manifolds using the Cheeger--Gromov invariant \cite{CG85} which associates to a
pair $(N,\pi_1(N)\to H)$, $H$ a group, a real number $\rho^{(2)}(N,\pi_1(N)\to H)$.  Using
standard arguments (cf. e.g. \cite{COT04}) one can show that
\[ \ba{rcl} &&\rho^{(2)}(N(K_i),\pi_1(N(K_i))\to \pi_1(N(K_i))/\pi_1(N(K_i))^{(2)})\\
&=& \rho^{(2)}(N(K),\pi_1(N(K))\to \pi_1(N(K))/\pi_1(N(K))^{(2)})+
\rho^{(2)}(N(C_i),\pi_1(N(C_i))\to \Z).\ea\] But $\rho^{(2)}(N(C_i),\pi_1(N(C_i))\to \Z)$ equals
the integral over the Levine--Tristram signature of $C_i$, which equals $\frac{4}{3}i$ (cf. e.g.
\cite[Lemma~5.3]{COT04}).

 \end{proof}

%===========================================
\section{The main theorem for finite fields} \label{section:finitefields}
As we have seen in Section \ref{sectionapp}, under favorable circumstances Proposition
\ref{sympcorg} can be applied directly to verify that for certain manifolds $N$, $S^{1} \times
N$ admits no symplectic structure. For other cases, however, there is no simple way to compute
in a compact way twisted Alexander polynomials, and we are led to \textit{brute force}
calculations, based on Fox Calculus  (see \cite{CF77}). The goal of this section is to determine
an effective method for implementing the results of Section \ref{fini} to get constraints that
can be explicitly verified by a computer algebra system. As our main concern at this point is
computational, we will focus only on the consequences of Theorem \ref{mainthm}, but it should
be clear that the results discussed in what follows can be formulated more or less
\textit{verbatim} for the study of fibered classes, partially recovering this way the results of
\cite{FK06}.

The one--variable Alexander polynomial $\Delta_{\phi}^{\a} \in \zt$  can in theory be computed
from a presentation of $\pi_1(N)$ using Fox calculus. In practice this turns out to be
non viable since this requires the computation of determinants of large matrices over the ring
$\zt$ and since it requires the computation of many greatest common divisors.

It is much easier to compute determinants and hence twisted Alexander polynomials over the
principal ideal domain $\F_{p}[t^{\pm 1}]$. As we will see switching to $\F_{p}[t^{\pm 1}]$ also
has the added advantage that we will be able to work with smaller representations (cf. Theorem
\ref{mainthmgen}).

Let $N$ be a 3--manifold, $\a:\pi_1(N) \to G$ an epimorphism to a finite group, and let $\phi \in
Hom(H,F)$ be a non--trivial homomorphism. We denote the induced
representation $\pi_1(N)\to \aut_\Z(\Z[G])$ (cf. Section \ref{twistal} for details) by $\a$ as
well, and given a prime $p$ we denote the representation $\pi_1(N)\to \aut_{\F_p}(\F_p[G])$ by
$\a_p$. We therefore get Alexander polynomials $\Delta^\a\in \Z[F]$ and $\Delta^{\a_p}\in
\F_p[F]$.
Remember that $\Delta^{\a_{p}} \in \F_{p}[F]$ is well-defined only up to an overall multiplication
by a non--zero element of $\F_{p}$, see Section \ref{twistal}. In particular we lose the notion
of monicness. This polynomial is related to $\Delta^{\a}$ by the following proposition, which is
another direct corollary to Lemma \ref{lem:alex02} and Theorem \ref{thm:functorial}.

\begin{proposition} \label{red}
Let $\phi:H\to F$ be a non--trivial  homomorphism to a free abelian group $F$. Let $\Delta^{\a} \in
\Z[F]$ and $\Delta^{\a_{p}} \in \F_{p}[F]$ be the twisted Alexander polynomials associated to an
epimorphism $\a: \pi_{1}(N) \rightarrow G$. Then
\begin{equation} \label{redp} \Delta^{\a_p} = \Delta^{\a} \ \ \mbox{mod $p$} \end{equation}
where the relation above holds in $\F_{p}[F]$ up to multiplication by a unit in $\F_{p}[F]$.
\end{proposition}

Proposition \ref{red} provides a way to determine, at least in principle, the twisted Alexander
norm $\| \cdot \|_{A}^{\a}$ in terms of the family of norms $\| \cdot \|_{A}^{\a_{p}}$
associated to the twisted Alexander polynomials $\Delta^{\a_{p}}$ since Equation (\ref{redp})
implies that $\mbox{supp} \Delta^{\a_{p}} \subset \mbox{supp}\Delta^{\a}$ and $\bigcup_{p}
\mbox{supp} \Delta^{\a_{p}} = \mbox{supp}\Delta^{\a}$. Furthermore Proposition \ref{red} allows
us  to characterize the elements of $\mbox{supp}\Delta^{\a}$ with coefficient $\pm 1$, as these
are the only ones whose reduction $mod$ $p$ is never trivial. In sum, we have the following
corollary:

\begin{corollary} \label{cofi}
Let $N$ be a closed 3--manifold and let $\a:\pi_1(N)\to G$ an epimorphism to a finite group.
Then for all $\phi \in H^{1}(N;\R)$
\[ \| \phi \|^{\a}_{A} = \max_{p}\{ \| \phi \|^{\a_p}_{A}\}. \]
Furthermore if $\phi$ lies in the cone over a top--dimensional face of the unit ball of the twisted
Alexander norm $||-||_A^\a$ which is dual to an element of $\mbox{supp}\Delta^{\a}$ with
coefficient $\pm 1$, then $\| \phi \|^{\a}_{A} = \| \phi \|^{\a_p}_{A}$ for all primes $p$.
\end{corollary}

\noindent Using Corollary \ref{cofi}, we can write a reformulation of Theorem \ref{strongsymp}
for finite fields:

\begin{theorem} \label{fifi} Let $N$ be an irreducible $3$--manifold such that $S^{1} \times N$
admits a symplectic structure with K\"unneth component $\varphi \in H^{1}(N)$. Then there exists a
top--dimensional face $F_{T}$ of the Thurston unit ball such that $\varphi \in \R_{+}F_{T}$.
Furthermore  for any prime $p$ we have $\Delta^{\a_p}\ne 0$ and for all $\phi
\in \R_{+}F_{T} \cap H^{1}(N)$ we have
\begin{equation} \label{normp} \| \phi \|^{\a_p}_{A}  = |G| \, \| \phi \|_{T} + \left\{
\begin{array}{ll} 0 & \mbox{if $b_{1}(N) >
     1$},  \\ \\ 2 \div\, \phi_{G}  & \mbox{if $b_{1}(N) =
     1$}. \end{array} \right. \end{equation} In particular, if $\phi \in \R_{+}F_{T} \cap H^{1}(N)$,
$\deg\, \Delta_{\phi}^{\a_p} = |G| \, \|\phi\|_{T} + 2 \div \, \phi_{G}$ for all primes
$p$.
\end{theorem}

\begin{remark}
\bn
\item Note that, due to the indeterminacy in the definition of the coefficients of
$\Delta^{\a_p}$, condition (1) in Theorem \ref{strongsymp} implies just that $\lambda_{G}$ has
non--zero coefficient in $\Delta^{\a_p}$. In this sense, that condition is automatically
included in Equation (\ref{normp}). On the other direction, the non--vanishing of the
coefficient of $\lambda_{G}$ in $\Delta^{\a_p}$ for all primes $p$ (implied by Equation
(\ref{normp})) guarantees that the coefficient of $\lambda_{G}$ in $\Delta^{\a}$ equals $\pm 1$.
\item Theorem \ref{fifi} ensures that if $N$ is an
irreducible manifold and $\phi \in \R_{+}F_{T} \cap H^{1}(N)$ represents the K\"unneth component
of a symplectic form, then $\Delta^{\a_p}_{\phi}$ is non--trivial for all primes $p$. \en
\end{remark}

At this point, we could implement Theorem \ref{fifi} by computing the polynomials
$\Delta^{\a_p}_{\phi}$ using Fox calculus. However, the computations would still be very time
consuming since the matrices involved are very large. It is therefore much more useful to find
criteria involving `small' representations of $\pi_1(N)$: Assume we are given a representation
$\b_p:G\to \aut_{\F_p}(B)$, $B$ a finite $\F_p$--module. We have Alexander polynomials
$\Delta_{N,\phi,i}^{\b_p\circ \a}\in \F_p\tpm$ which we denote by $\Delta_{N,\phi,i}^{\b_p}$. It
is well--known that given a presentation of $\pi_1(N)$ these polynomials can be computed  for
$i=0,1$ using Fox calculus; we refer to \cite{FK06} for the computation when $i=2$.

For these Alexander polynomials, a
McMullen-type inequality (similar to Lemma \ref{normineq}) is one of the main results of
\cite{FK05} and \cite{FK06}:

\begin{theorem} \label{mainthmfk05} Let $N$ be a  3--manifold, $B$ a finite  $\F_p$--module and let $\b_p:G\to
\mbox{Aut}_{\F_p}(B)$ be a representation. If $\Delta^{\b_p}\ne 0$, then for all non--zero $\phi
\in H^1(N)$ we have the following inequalities:
\[ \| \phi \|^{\b_p}_{A}  \leq \mbox{dim}(B) \| \phi \|_{T} + \left\{ \begin{array}{ll} 0 & \mbox{if $b_{1}(N) >
     1$},  \\ \\  \deg \Delta^{\b_p}_{\phi,0} + \deg \Delta^{\b_p}_{\phi,2}  & \mbox{if $b_{1}(N) =
     1$}. \end{array} \right. \]
     \end{theorem}
Note that in the case where $B = \F_{p}[G]$ we can apply Lemma \ref{lem:alex02} to conclude that
$\div \phi_{G} = \deg \Delta^{\a_p}_{\phi,0} = \deg \Delta^{\a_p}_{\phi,2}$.

The constraints on the existence of symplectic structures obtained thus far (Theorems
\ref{strongsymp} and \ref{fifi}) apply to the twisted Alexander polynomials defined using the
modules $\Z[G][H]$ and $\F_{p}[G][H]$. We want to translate them into conditions on the
polynomials $\Delta_{N,\phi,i}^{\b_p}$. To do so, we will need a further assumption, namely we
will choose $\F_{p}$ and $G$ so that any $\F_{p}[G]$--module is semisimple, i.e. that the
$\F_{p}$--representations of $G$ are \textit{completely reducible}. Concerning this we have the
following well--known version of Maschke's Theorem:

\begin{theorem}\label{thmmaschke}
Let $p$ be a prime such that  $(p,|G|) = 1$. Then every $\F_p[G]$--module is semisimple.
Moreover, every simple $\F_p[G]$--module is a direct summand of $\F_p[G]$.
\end{theorem}

\begin{proof}
For the first part we refer to \cite[p.~216]{Ro96}. For the second part, let $C$ be a simple
$\F_p[G]$--module and let $v\in C$. Then there exists a short exact sequence of
$\F_p[G]$--modules:
\[ \ba{ccccccccc} 0 &\to &K&\to& \F_p[G]&\to & C&\to&0\\
&&&&g&\mapsto&vg.&\ea \]
Note that the last map must be an epimorphism since $C$ is simple.
Since $\F_p[G]$ is a semisimple left $\F_p[G]$--module, the sequence above splits and  we can
write $\F_p[G]=K \oplus C$.
\end{proof}

Let us then assume that $(p,|G|) = 1$. The twisted Alexander polynomials associated to a
reducible $\F_{p}[G]$--module satisfy the following lemma.
\begin{lemma} \label{decompose}
Let $\b_p:G\to \aut_{\F_p}(B)$ and $\g_p:G\to \aut_{\F_p}(C)$ be finite dimensional
representations. Denote the direct sum of the representations by $\b_p\oplus \g_p:G\to
\aut_{\F_p}(B\oplus C)$. Furthermore let $\phi:\pi_1(N)\to F$ be a non--trivial homomorphism to
a free abelian group $F$. Then
\[ \Delta^{\b_p \oplus \g_p}_{\phi,i} = \Delta^{\b_p}_{\phi,i} \cdot
\Delta^{\g_p}_{\phi,i}\] and, for all $\phi \in H^{1}(N;\R)$, $\|\phi \|_{A}^{\b_p \oplus \g_p}
= \|\phi \|_{A}^{\b_p} + \|\phi \|_{A}^{\g_p}$.
\end{lemma}

\begin{proof}
We have $H_*(N;(B\oplus C)[F]) \cong H_*(N;B[F]) \oplus H_*(N;C[F])$; therefore we can get a free
resolution of $H_*(N;(B\oplus C)[F])$ as direct sum of the resolutions of $H_*(N;B[F])$ and
$H_*(N;C[F])$, which implies the relation between the twisted Alexander polynomials. The
relation between the twisted Alexander norms follows at once from that relation (for the case of
$F = H$), for example by looking at the dense subset of $H^{1}(N;\R)$  given by the intersection
of the cones on the top--dimensional faces of the unit balls of the norms $\| \cdot
\|_{A}^{\b_p}$ and  $\|\cdot \|_{A}^{\g_p}$, and using continuity. \end{proof}

\noindent We can state now the analogue of Theorem \ref{fifi} for a general $\F_{p}[G]$--module
$B$ with $(p,|G|)=1$:

\begin{theorem}\label{mainthmgen}
Let $N$ be an irreducible $3$--manifold such that $S^{1} \times N$ admits a symplectic structure
with K\"unneth component $\varphi \in H^{1}(N)$; then there exists a top--dimensional face
$F_{T}$ of the Thurston unit ball such that $\varphi \in \R_{+}F_{T}$ and such that for all
$\phi \in \R_{+}F_{T} \cap H^{1}(N)$, any epimorphism $\a:\pi_1(N)\to G$ onto a finite group $G$ and any $\F_{p}[G]$--module $B$ with $(p,|G|)=1$
\begin{equation}
\label{mainfor} \| \phi \|^{\b_p}_{A}  = \mbox{dim}(B) \| \phi \|_{T} + \left\{
\begin{array}{ll} 0 & \mbox{if $b_{1}(N) >
     1$},  \\ \\ \deg \, \Delta^{\b_p}_{\phi,0} + \deg \, \Delta^{\b_p}_{\phi,2} & \mbox{if $b_{1}(N) =
     1$}.  \end{array} \right. \end{equation} In particular, if $\phi \in \R_{+}F_{T} \cap H^{1}(N)$,
$\deg\, \Delta_{\phi}^{\b_p} = \mbox{dim}(B) \|\phi\|_{T} + \deg \, \Delta^{\b_p}_{\phi,0}
+ \deg \,\Delta^{\b_p}_{\phi,2} $ for all primes $p$. \end{theorem}

\begin{proof}
First observe that by Maschke's Theorem and Lemma \ref{decompose}, we can limit ourselves to
consider the case where $B$ is a simple $\F_{p}[G]$--module, and thus a direct summand of
$\F_{p}[G]$. Denote by $C$ the $\F_p[G]$--module such that $\F_{p}[G] = B \oplus C$.  Denote the
$\F_p[G]$--module structures by $\b_p$ respectively $\g_p$. Note that it follows from Theorem
\ref{fifi} and Lemma \ref{decompose} that $\Delta^{\b_p}\ne 0$ and $\Delta^{\g_p}\ne 0$, which
allows us to apply Theorem \ref{mainthmfk05}. Now applying in order
Theorem \ref{fifi}, Lemma \ref{decompose} and Theorem \ref{mainthmfk05} (twice) we have
\[ \ba{rcl} &\hspace{-0.15cm}&\hspace{-0.15cm}|G| \, \| \phi \|_{T} + \left\{ \begin{array}{l} 0  \\ \\ \deg \Delta^{\a_p}_{\phi,0} + \deg \Delta^{\a_p}_{\phi,2} \end{array} \right. = \|\phi \|_{A}^{\a_p} = \|\phi \|_{A}^{\b_p} + \|\phi \|_{A}^{\g_p} \leq
\\ \\ &\hspace{-0.15cm}\leq&\hspace{-0.15cm} ||\phi||_A^{\b_p}+ \mbox{dim}(C) \| \phi \|_{T} + \left\{ \begin{array}{ll} 0 &\hspace{-0.15cm} \mbox{if $b_{1}(N) >
     1$},  \\ \\  \deg \Delta^{\g_p}_{\phi,0} + \deg \Delta^{\g_p}_{\phi,2} &\hspace{-0.15cm} \mbox{if $b_{1}(N) =
     1$};  \end{array} \right.
\\ \\ &\hspace{-0.15cm}\leq&\hspace{-0.15cm} (\mbox{dim}(B) + \mbox{dim}(C)) \| \phi \|_{T} + \left\{ \begin{array}{ll} 0 &\hspace{-0.15cm} \mbox{if $b_{1}(N) >
     1$},  \\ \\ \deg \Delta^{\b_p}_{\phi,0} + \deg \Delta^{\b_p}_{\phi,2} + \deg \Delta^{\g_p}_{\phi,0} + \deg \Delta^{\g_p}_{\phi,2} &\hspace{-0.15cm} \mbox{if $b_{1}(N) =
     1$}.  \end{array} \right. \ea \]
It follows from Lemma \ref{decompose} that the last term equals the first term.  But this implies
that all inequalities are in fact equalities, and this shows that (\ref{mainfor}) holds for a simple
module, and thus for all $\F_{p}[G]$--modules. Finally, the condition on $\deg\Delta_{\phi}^{\b_p}$
follows from (\ref{mainfor}), a version of Theorem
\ref{thm:functorial} for general representations proved in \cite{FK05} and the argument in the proof of Proposition
\ref{holdt}, applied to $\Delta^{\b_p}$.
\end{proof}

\begin{remark} In studying the problem of fiberability of $(N,\phi)$ we could have obtained the result of Theorem \ref{mainthmgen} directly for \textit{any} $R[G]$--module, without having to pass through the results of Section \ref{fini} (this is in fact the approach of \cite{FK06}). Instead, for the symplectic case, we have to follow that path, as the twisted Alexander polynomial $\Delta^{\a}$ associated to $R[G]$ is the only one that can be interpreted directly in terms of Seiberg-Witten theory, allowing us to apply the corresponding symplectic constraints. It is entirely possible, however, that Theorem \ref{mainthmgen} holds true without the condition $(p,|G|)=1$. \end{remark}
%===========================================
\section{Kronheimer's example and 12--crossing knots} \label{section:ex}

%============================
\subsection{Kronheimer's example: the  Pretzel knot $P=(5,-3,5)$}
\begin{figure}[h] \centering
\includegraphics[scale=0.25]{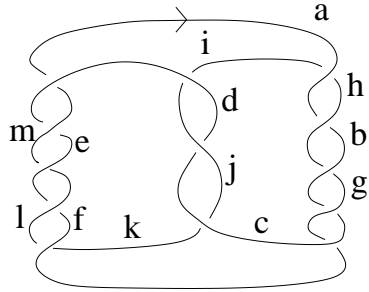}
\caption{The pretzel knot $(5,-3,5)$.} \label{pretzel} \end{figure}

We will now consider the Pretzel knot $P=(5,-3,5)$, whose diagram is given in Figure \ref{pretzel}.
This knot has Alexander polynomial $t^2-3t+1$ and has a Seifert surface of genus one.
 The fundamental group of
$S^3\sm \nu P$ is generated by the meridians $a,b,\dots,m$ of the strands in the knot diagram of
Figure \ref{pretzel}. The relations are
\[ \ba{rclrclrclrclrcl}
a&=&h^{-1}bh,&    b&=&g^{-1}cg,&    c&=&jdj^{-1},&    d&=&m^{-1}em,&    e&=&l^{-1}fl,\\
f&=&c^{-1}gc,&    g&=&b^{-1}hb,&    h&=&a^{-1}ia,&    i&=&djd^{-1},&    j&=&ckc^{-1},\\
k&=&f^{-1}lf,&    l&=&e^{-1}me,&    m&=&d^{-1}ad.& \ea \]
The longitude is represented by
$\l=h^{-1}g^{-1}jm^{-1}l^{-1}c^{-1}b^{-1}a^{-1}dcf^{-1}e^{-1}d^{-1}a^7$. Clearly
$\pi_1(N(P))=\pi_1(S^3\sm \nu P)/\langle \l \rangle$.

Using the computer program \emph{KnotTwister} \cite{Fr05} we found an epimorphism  $\a:\pi_1(N(P))\to  S_5$ which is given by
\[ \begin{array}{rclrclrclrcl}
a &\mapsto& (51234)&b &\mapsto& (43521)&c &\mapsto& (54132)&d &\mapsto& (43521)\\
e &\mapsto& (35214)&f &\mapsto& (23451)&g &\mapsto& (35214)&h &\mapsto& (24153)\\
i &\mapsto& (35421)&j &\mapsto& (41532)&k &\mapsto& (54213)&l &\mapsto& (41532)\\
m &\mapsto& (24153).
\end{array} \]
Now consider the representation $\b_{7}:\pi_1(N(P))\xrightarrow{\a}
S_5\to \gl(\F_{7},5)$ where $S_5$ acts on $\F_{7}^5$ by permutation of the coordinates. Using
\emph{KnotTwister} we compute that $\Delta^{\b_7}=0$ and $\deg \Delta_0^{\b_7}=\deg
\Delta_2^{\b_7}=1$.  Since $p=7$ and $|S_5|=120$ are coprime it follows immediately from Theorem
\ref{mainthmgen}, applied to a generator $\phi \in H^{1}(N(P))$, that $S^1\times N(P)$ is not
symplectic.

%============================
\subsection{12--crossing knots} It is well--known that a knot $K$ with 11 crossings or less is
fibered (and equivalently $N(K)$ is fibered) if and only if $2\, g(K)=\deg(\Delta_K)$ and
$\Delta_K$ is monic. It follows therefore from Theorem \ref{thmvidussi} that if $K$ is a knot
with 11 crossings or less, then $S^1 \times N(K)$ is symplectic if and only if $K$ is fibered.

By work of Stoimenow and Hirasawa and the results of   \cite{FK06}  there exist 13 non--fibered
12--crossing knots which have monic Alexander polynomials and such that $\deg(\Delta_K)=2\,
g(K)$. According to \cite{FK06} these knots are $12_{1345}$, $12_{1498}$, $12_{1502}$,
$12_{1546}$, $12_{1567}$, $12_{1670}$, $12_{1682}$, $12_{1752}$, $12_{1771}$, $12_{1823}$,
$12_{1938}$, $12_{2089}$, $12_{2103}$ (where we use \emph{Knotscape} notation \cite{HT}). For
these 13 knots Theorem \ref{thmvidussi} is not strong enough to show that $S^1\times N(K)$ is
not symplectic. The following table shows all the knots for which we found a representation
$\b_p:\pi_1(N(K))\to G=S_k\to \gl(\F_p,k)$ (where $S_k$ acts by permutation on $\F_p^k$) such
that $\gcd(p,|S_k|)=1$ and such that either $\Delta^{\b_p}=0$ or
\[ \deg \Delta^{\b_p} < k \cdot (2g(K) -2) + \deg \Delta^{\b_p}_0 + \deg \Delta^{\b_p}_2. \]
For these knots $K$ the manifold $S^1\times N(K)$ is not symplectic by Theorem \ref{mainthmgen}.
In the following table we give all the relevant data.
Recall that we write $\deg \Delta^{\b_p}=-\infty$ if
$\Delta^{\b_p}=0$.
\[
\ba{|r|r|r|r|r|r|r|} \hline \mbox{Knot} &12_{1345}& 12_{1498}& 12_{1502}& 12_{1546}& 12_{1567}&
12_{1752} \\ \hline
\mbox{$\| \phi \|_{T} = 2 g(K) -2$} &2&4&4&2&2&2\\
\mbox{Order $k$ of permutation group}&6&6&5&5&5&6\\
\mbox{Prime $p$ }&7&7&11&7&7&17\\
\deg \Delta^{\b_p}_0 =\deg \Delta^{\b_p}_2 &1&1&1&1&1&1\\
\deg \Delta^{\b_p} &-\infty&24&14&-\infty&-\infty&10\\
k(2g(K)-2)+\deg \Delta_0^{\b_p}+\deg \Delta_2^{\b_p}& 14&26&22&12&12&14\\
\hline \mbox{Knot} &12_{1670}& 12_{1771}&\ 12_{1823}& 12_{1938}& 12_{2089}& 12_{2103}\\ \hline
\mbox{$\| \phi \|_{T} = 2 g(K) -2$} &2&2&2&2&2&2\\
\mbox{Order $k$ of permutation group}&6&5&6&5&5&5\\
\mbox{Prime $p$ }&17&7&7&11&11&7\\
\deg \Delta^{\b_p}_0=\deg \Delta_2^{\b_p}  &1&2&2&1&1&1\\
\deg \Delta^{\b_p} &10&10&-\infty&4&4&-\infty\\
k(2g(K)-2)+\deg \Delta_0^{\b_p}+\deg \Delta_2^{\b_p}& 14&16&16&14&14&14\\
\hline \ea \]

This approach did not work for the knot  $K = 12_{1682}$. The stumbling block in our efforts was
the condition $(|G|,p)=1$. However, we can bypass this problem by applying directly Theorem
\ref{fifi} using a map $\a:\pi_1(N(K))\to A_4$ such that $\Delta^{\a_3}\in \F_{3}[t^{\pm 1}]$
has degree 21 and $\deg(\Delta_{0}^{\a_3})=\deg(\Delta_{2}^{\a_3})=3$. Note that $|A_4|=12$ and
that $g(K)=2$. It follows that
\[ \deg \Delta^{\a_3}=21 \ne 30=|A_4| (2g(K)-2)+\deg \Delta_{0}^{\a_3}+\deg\Delta_{2}^{\a_3}.\]
Therefore $S^1\times N(K)$ is not symplectic by Theorem \ref{fifi}. Together with the previous
results, we conclude that Conjecture \ref{conjfolk} holds for all knots with up to $12$ crossings.

\begin{remark}
Note that the proof of non--fiberedness of these knots contained in \cite{FK06} does not, in
many cases, apply to show that $S^{1} \times N(K)$ is not symplectic, as the (computationally
somewhat simpler) representations therein considered mostly fail to satisfy the condition
$(p,|G|)=1$.
\end{remark}

%===========================================

%===========================================
\section{The fiberedness conjecture and 3--manifold groups}\label{sectionconj} \label{conjfibered}

The aim of this section is to discuss some evidence to Conjecture \ref{conjcha}. To start, it is
more or less explicitly known that the hypothesis of Conjecture \ref{conjcha} are sufficient
(assuming the Geometrization Conjecture) to show that $N$ is prime, a condition that is
well--known to hold for fibered manifolds. Precisely, we have the following:
\begin{proposition} \label{lem:prime}
Let $N$ be a 3--manifold which is not prime and let $\phi\in H^1(N)$. If the Geometrization
Conjecture holds, then there exists an epimorphism $\a: \pi_1(N)\to G$, $G$ a finite group, such
that $\twialexg=0\in \zt$.
\end{proposition}
This statement affords a direct proof using a Mayer--Vietoris decomposition for the Alexander
module for a connected sum. It can also be simply deduced (see \cite{McC01} or \cite{Vi99}) from
standard vanishing theorems for Seiberg-Witten invariants of a suitable cover of $N$, whose
existence is granted by the Geometrization Conjecture, which implies that $\pi_1(N)$ is
residually finite (cf. \cite{He87}). Note that without that assumption, the proof shows that a
manifold $N$ satisfying the hypothesis of Conjecture \ref{conjcha} is prime up to the connected
sum with a homology sphere which does not admit non--trivial finite covers.

In the following let $N$ be a 3--manifold and let $\phi \in H^{1}(N)$ a primitive class such
that for \textit{any} epimorphism onto a finite group $\a: \pi_{1}(N) \rightarrow G$ the twisted
Alexander polynomial $\Delta^{\a}_{N,\phi} \in \Z[t^{\pm 1}]$ is monic and
$\deg\Delta_{N,\phi}^{\a} = |G| \, \|\phi\|_{T} + 2 \div \phi_{G}$. Let $S$ be a Thurston
norm minimizing surface dual to $\phi$. As by hypothesis $\Delta_{\phi} \neq 0$ we can assume,
by \cite{McM02}, that $S$ is connected. Then let $M$ be the result of cutting $N$ along $S$.
Denote the positive and the negative inclusions of $S$ into $M$ by $i_+$ and $i_-$. Since $S$ is
minimal, $i_+: \pi_1(S)\to \pi_1(M)$  is injective by Dehn's Lemma. Since $N$ is prime it
follows easily from Stallings' theorem \cite{St62} that $(N,\phi)$ fibers over $S^1$ if and only
if $i_+: \pi_1(S)\to \pi_1(M)$ is surjective.

\begin{proposition}\label{propdeltamodp}
$H_1^\a(S;\Z[G])$ and $H_1^\a(M;\Z[G])$ are free abelian groups of the same rank.
\end{proposition}

\begin{proof}
First note that $H_1^\a(S;\Z[G])$ is the first homology of the (not necessarily connected) cover
of $S$ corresponding to $\pi_1(S)\to G$, hence $H_1^\a(S;\Z[G])$ is free abelian.

Let $p$ be a prime. It follows from Proposition \ref{red} that $\Delta^{\a_p}_{\phi}\ne 0$. By
Lemma \ref{lem:alex02} $H_2^\a(N;\fp[G][t^{\pm 1}])$ is $\fpt$--torsion. Therefore the
Mayer--Vietoris sequence for $N$ with $\fp[G]\tpm$--coefficients gives a short exact sequence
(cf. \cite{FK06} for details)
\[ 0\to H_1^{\a_p}(S;\fp[G])\otimes_{\F_p} \fpt \xrightarrow{ti_+-i_-}  H_1^{\a_p}(M;\fp[G])\otimes_{\F_p}
\fpt \to  H_1^{\a_p}(N;\fp[G][t^{\pm 1}])\to 0.\] Tensoring with $\F_p(t)$ we see that in
particular $H_1^{\a_p}(S;\fp[G])\cong H_1^{\a_p}(M;\fp[G])$ for every prime $p$ as
$\F_p$--vector spaces. Note that $H_0^\a(S;\Z[G])\cong \Z[G/\im\{\pi_1(S)\to G\}]$ and
$H_0^\a(M;\Z[G])\cong \Z[G/\im\{\pi_1(M)\to G\}]$; in particular both are $\Z$--torsion free. It
follows from the universal coefficient theorem applied to the $\Z$--module complex
$C_*(\ti{S})\otimes_{\Z[\pi_1(S)]}{\Z[G]}$ that
\[ H_1^\a(S;\Z[G])\otimes_\Z \F_p \cong H_1^{\a_p}(S;\F_p[G])\]
for every prime $p$. The same statement holds for $M$. Combining our results we see that
$H_1^\a(S;\Z[G])\otimes_\Z \F_p\cong H_1^\a(M;\Z[G])\otimes_\Z \F_p$ for any prime $p$. It
follows that $H_1^\a(S;\Z[G])\cong H_1^\a(M;\Z[G])$.
\end{proof}

\noindent Now consider the exact sequence
\[ H_1^\a(S;\Z[G])\otimes \zt \xrightarrow{ti_+-i_-} H_1^\a(M;\Z[G])\otimes \zt \to H_1^\a(N;\Z[G]\tpm)\to 0.\]
Since $H_1^\a(S;\Z[G])$ and $H_1^\a(M;\Z[G])$ are free abelian groups of the same rank it
follows that $\twialexg=\det(ti_+-i_-)$. Note that this does not yet imply that $i_+$ and $i_-$
have full rank. Write $b_i^\a(S)=\rank_{\Z}(H_i^\a(S;\Z[G]))$. A standard Euler characteristic
argument shows that
\[-b_0^\a(S)+b_1^\a(S)-b_2^\a(S)=-|G|\chi(S)=|G|\cdot ||\phi||_T.\]
Furthermore $b_i^\a(S)=\deg \Delta^{\a}_{\phi,i}$ for $i=1,2$ by
\cite[Propositions~4.9~and~4.11]{FK06}. Combining this with  the assumption
\[ \deg\Delta_{\phi}^{\a} = |G| \, \|\phi\|_{T} + 2 \div \phi_{G}=|G| \, ||\phi||_T+
\deg\Delta_{\phi,0}^{\a}+\deg\Delta_{\phi,2}^{\a}\] we get that
$\deg(\det(ti_+-i_-))=\deg \Delta^{\a}_{\phi}=b_1^\a(S)$. Since $i_+$ and $i_-$ are
$b_1^\a(S)\times b_1^\a(S)$ matrices over $\Z$ it now follows clearly that $\det(i_+)$ equals the top
coefficient of $\Delta^\a_{\phi}$, namely $1$. Similarly $\det(i_-)=1$. This shows that $i_+:
H_1^\a(S;\Z[G])\to H_1^\a(M;\Z[G])$ is an isomorphism.

Summarizing, this discussion shows that Conjecture \ref{conjcha} follows from the Geometrization
Conjecture and the following question in the theory of 3--manifold groups.

\begin{conjecture} \label{conjcha2} Let $S$ be an incompressible surface in a 3--manifold $N$
and let $M$ be $N$ cut along $S$. Let $i:S\to M$ be one of the positive and the negative
inclusions of $S$ into $M$. If $i: H_1^\a(S;\Z[G])\to H_1^\a(M;\Z[G])$ is an isomorphism for
every homomorphism $\pi_1(N)\to G$, $G$ a finite group, then $i: \pi_1(S)\to \pi_1(M)$ is
surjective.
\end{conjecture}

Note that it is well--known that the inclusion induced homomorphisms $\pi_1(S) \to \pi_1(N)$ and
$\pi_1(M) \to \pi_1(N)$ are injections. We think that an affirmative answer to the above
conjecture will need a strong condition on $\pi_1(N)$ (of which $\pi_1(M)$ and $\pi_1(S)$ are
subgroups) such as \emph{subgroup separability}, which is conjectured to hold for fundamental
groups of hyperbolic manifolds (cf. \cite{Th82}).


\begin{thebibliography}{10}
\bibitem[Ba01]{Ba01} S. Baldridge, {\em  Seiberg-Witten invariants of 4-manifolds with free circle actions},
 Commun. Contemp. Math.  3: 341--353 (2001).
\bibitem[Ch03]{Ch03} J. Cha, {\em Fibred knots and twisted Alexander invariants},
Transactions of the AMS 355: 4187--4200 (2003)
\bibitem[CG85]{CG85}
J. Cheeger, M. Gromov, {\em Bounds on the von Neumann dimension of $L\sp 2$-cohomology and the
Gauss-Bonnet theorem for open manifolds},  J. Differential Geom. 21,  no. 1: 1--34 (1985)
\bibitem[COT04]{COT04}
T. Cochran, K. Orr, P. Teichner, {\em Structure in the classical knot concordance group},
Comment. Math.
Helv.  79: 105-123 (2004)
\bibitem[CF77]{CF77} R. Crowell,
R. Fox, {\em Introduction to knot theory}, Graduate Text in Mathematics 57, Springer Verlag
(1977)
\bibitem[Du01]{Du01} N. Dunfield, {\em Alexander and Thurston norm of fibered 3--manifolds}, Pacific J. Math.,
200, no. 1: 43--58 (2001)
\bibitem[Du85]{Du85}
M. J. Dunwoody, {\em An equivariant sphere theorem}, Bull. London Math. Soc. 17, no. 5, 437--448
(1985)
\bibitem[Fr05]{Fr05}
S. Friedl, {\em KnotTwister}, \texttt{http://www.labmath.uqam.ca/\~\,friedl/index.html} (2005).
\bibitem[FK05]{FK05} S. Friedl, T. Kim, {\em  Twisted Alexander norms give lower
bounds on the Thurston norm}, Trans. Amer. Math. Soc. (to appear)
\bibitem[FK06]{FK06} S. Friedl, T. Kim,
{\em The Thurston norm, fibered manifolds and twisted Alexander polynomials}, Topology 45, 929--953 (2006)
\bibitem[FV07]{FV07}
S. Friedl, S. Vidussi, {\em Symplectic $\bf S^{1} \times N^3$,
surface subgroup separability, and totally degenerate Thurston
norm}, Preprint (2007), to appear in J. Amer. Math. Soc.
\bibitem[Ga83]{Ga83} D. Gabai, {\em Foliations and the
topology of 3--manifolds}, J. Differential Geometry 18, no. 3: 445--503 (1983)
\bibitem[Ga87]{Ga87b} D. Gabai, {\em Foliations
and the topology of 3--manifolds. III}, J. Differential Geometry 26, no. 3: 479--536 (1987)
\bibitem[GKM05]{GKM05} H. Goda, T. Kitano, T. Morifuji, {\em Reidemeister Torsion, Twisted Alexander Polynomial and Fibred Knots}, Comment. Math. Helv.  80,  no. 1: 51--61 (2005)
\bibitem[GS99]{GS99} R. Gompf, A. Stipsicz, {\em 4-manifolds and Kirby calculus}, Graduate Studies
in Mathematics 20, AMS (1999)
\bibitem[He87]{He87} J. Hempel, {\em Residual
finiteness for $3$-manifolds}, Combinatorial group theory and topology (Alta, Utah, 1984),
379--396, Ann. of Math. Stud., 111, Princeton Univ. Press, Princeton, NJ (1987)
\bibitem[HT]{HT}
J. Hoste, M. Thistlethwaite, {\em Knotscape},
\texttt{http://www.math.utk.edu/\~\,morwen/knotscape.html}
\bibitem[JW93]{JW93} B. Jiang, S. Wang, {\em Twisted topological invariants associated with representations}, Topics in
knot theory (Erzurum 1992), 211--227, NATO Adv. Sci. Inst. C. Math. Phys. Sci, 399, Kluwer Acad.
Publ., Dordrecht (1993)
\bibitem[KL99]{KL99a} P. Kirk, C. Livingston, {\em Twisted Alexander invariants, Reidemeister
torsion and Casson--Gordon invariants}, Topology 38, no. 3: 635--661 (1999)
\bibitem[Kr98]{Kr98} P. Kronheimer, {\em Embedded surfaces and gauge theory in three
and four dimensions},  Surveys in differential geometry, Vol. III (Cambridge, MA, 1996),
243--298, Int. Press, Boston, MA (1998)
\bibitem[Kr99]{Kr99} P. Kronheimer, {\em Minimal genus in $S\sp 1\times
M\sp 3$},  Invent. Math.  135,  no. 1: 45--61 (1999)
\bibitem[Le67]{Le67} J. Levine, {\em A method for generating link polynomials}, Amer. J. Math. 89
(1967) 69--84.
\bibitem[L01]{L01} X. S. Lin, {\em Representations of knot groups and twisted
Alexander polynomials}, Acta Math. Sin. (Engl. Ser.)  17,  no. 3: 361--380 (2001)
\bibitem[McC01]{McC01} J. McCarthy, {\em On the asphericity of a symplectic $M^{3} \times S^{1}$}, Proc. Amer. Math. Soc. 129: 257--264 (2001)
\bibitem[McM02]{McM02} C. T. McMullen, {\em The Alexander polynomial of a 3--manifold and the Thurston
norm on cohomology}, Ann. Sci. Ecole Norm. Sup. (4) 35, no. 2: 153--171 (2002)
\bibitem[MT99]{MT99} C. T. McMullen, C. H. Taubes, {\em 4-manifolds with inequivalent symplectic
forms and 3-manifolds with inequivalent fibrations}, Math. Res. Lett. 6: 681--696 (1999)
\bibitem[MeT96]{MeTa} G. Meng, C. H. Taubes, {\em SW = Milnor torsion}, Math. Res. Lett. 3: 661--674 (1996)
\bibitem[Ro96]{Ro96} D. Robinson, {\em A Course in the Theory of Groups}, Graduate
Texts in Mathematics 80, Springer Verlag (1996)
\bibitem[St62]{St62}
J. Stallings, {\em On fibering certain 3--manifolds}, 1962 Topology of 3--manifolds and related
topics (Proc. The Univ. of Georgia Institute, 1961) pp. 95--100 Prentice-Hall, Englewood Cliffs,
N.J. (1962)
\bibitem[Sto]{Sto}
A. Stoimenow, \texttt{http://www.ms.u-tokyo.ac.jp/\~\,stoimeno/ptab/index.html}
\bibitem[Ta98]{Ta98} C. H. Taubes, {\em The geometry of the Seiberg-Witten invariants}, Doc. Math. J. DMV Extra Volume ICM II: 493--492 (1998)
\bibitem[Th76]{Th76} W. P. Thurston,
{\em Some simple examples of symplectic manifolds}, Proc. Amer. Math. Soc. 55 (1976), no. 2, 467--468.
\bibitem[Th82]{Th82} W. P. Thurston,
{\em Three dimensional manifolds, Kleinian groups and hyperbolic geometry}, Bull. Amer. Math.
Soc. 6 (1982)
\bibitem[Th86]{Th86} W. P. Thurston, {\em A norm for the homology of 3--manifolds}, Mem.
Amer. Math. Soc. 339: 99--130 (1986)
\bibitem[Tu86]{Tu86} V. Turaev, {\em Reidemeister torsion in knot theory}, Russian Math. Surveys 41, no. 1: 119--182 (1986)
\bibitem[Tu01]{Tu01} V. Turaev, {\em Introduction to combinatorial torsions}, Birkh\"auser, Basel, (2001)

\bibitem[Vi99]{Vi99} S. Vidussi, {\em The Alexander norm is smaller than the Thurston norm; a
Seiberg--Witten proof}, Prepublication Ecole Polytechnique 6 (1999)
\bibitem[Vi03]{Vi03} S. Vidussi,
{\em Norms on the cohomology of a 3-manifold and SW theory},
 Pacific J. Math.  208,  no. 1: 169--186 (2003)
\bibitem[Wa94]{Wa94}
M. Wada, {\em Twisted Alexander polynomial for finitely presentable groups}, Topology 33, no. 2:
241--256 (1994)
\end{thebibliography}
\end{document}